\input amstex
\documentstyle{amsppt}
\NoBlackBoxes
\nologo \magnification1200 \pagewidth{6.5 true in}
\pageheight{9 true in} \topmatter
\title Distribution of values of $L$-functions at the edge of the critical strip
\endtitle
\author{ Youness Lamzouri }
\endauthor
\address D{\'e}partment  de Math{\'e}matiques et Statistique,
Universit{\'e} de Montr{\'e}al, CP 6128 succ Centre-Ville,
Montr{\'e}al, QC  H3C 3J7, Canada
\endaddress
\email{Lamzouri{\@}dms.umontreal.ca}
\endemail
\thanks
\endthanks
\abstract We prove several results on the distribution of values of $L$-functions at the edge of the critical strip, by constructing and studying a large class of random Euler products. Among new applications, we study families of symmetric power $L$-functions of holomorphic cusp forms in the level aspect (assuming the automorphy of these $L$-functions) at $s=1$,  functions in the Selberg class (in the height aspect),  and
quadratic twists of a fixed $GL(m)/{\Bbb Q}$-automorphic cusp form at $s=1$.
\endabstract
\endtopmatter
\define \ex{\Bbb E}
\define \pr{\text{Prob}}

\define \sums{\sideset \and^{\flat}\to\sum}
\define \sumh{\sideset \and^{h}\to\sum}
\document
\head {Introduction and statement of results }
\endhead

 \noindent  Values of $L$-functions at the edge of the critical strip are interesting objects that encode deep arithmetic information. For example, the non-vanishing of the Riemann zeta function $\zeta(s)$ on the line Re$(s)=1$ implies the prime number theorem proved by Hadamard and de La Vall\'ee Poussin in 1896. Furthermore, let $S_k^p(N)$ be the set of arithmetically normalized primitive cusp forms of weight $k$ and level $N$. Serre [27], proved that the Sato-Tate conjecture is equivalent to the non-vanishing of symmetric power $L$-functions $L(s,\text{Sym}^kf)$ at Re$(s)=1$, for all $k\in\Bbb{N}.$ This was recently used by Taylor {\it and al} [30] to prove the Sato-Tate conjecture for all non-CM elliptic curves over totally real fields, satisfying the mild condition of having multiplicative reduction at some prime.

 The distribution of these values have been extensively
studied over the last decades. One can quote the work of
Granville-Soundararajan [7] and Lamzouri [15] in the case of $\zeta(1+it)$; Elliott ([5] and [6]), Montgomery-Vaughan [20] and
Granville-Soundararajan [8]
 in the case of Dirichlet $L$-functions
of quadratic characters $L(1,\chi_d)$; Duke [4] in the case of Artin $L$-functions, and the work of
Cogdell-Michel [1], Habsieger-Royer [9], Lau-Wu [16], Liu-Royer-Wu [19], Royer ([21] and [22]), and Royer-Wu ([23] and [24]) in the case of symmetric power
$L$-functions of $GL(2)$-automorphic forms.

In the case of quadratic characters $L$-functions, Montgomery and Vaughan [20] were the first to construct a probabilistic model to study the distribution of this family at $s=1$, and based on their model they conjectured that the distribution function should be double exponentially decreasing. Granville and Soundararajan [8] computed large complex moments of $L(1,\chi_d)$ and deduced an asymptotic formula for the distribution of the $L$-values, proving Montgomery and Vaughan's conjecture. More precisely let $\Phi_x(\tau)$ to be the proportion of fundamental discriminants $d$ with $|d|\leq x$, for which $L(1,\chi_d)>e^{\gamma}\tau$. Then Granville and Soundararajan proved that uniformly in the region $\tau\leq \log_2x$ (throughout this paper
$\log_j$ denotes the $j$-th iterated logarithm, so that $\log_1n
=\log n$ and $\log_j n=\log (\log_{j-1}n) $ for each $j\geq 2$), we have
$$ \Phi_x(\tau)
=\exp\left(-\frac{e^{\tau-C_1}}{\tau}
\left(1+O\left(\frac{1}{\tau}\right)\right)\right),\tag{1}$$
where $C_1:=1+\displaystyle{\int_0^1\log\cosh y\frac{dy}{y^2} + \int_1^{\infty}(\log\cosh y-y)\frac{dy}{y^2}}=0.8187...$

Using a similar approach, they studied the distribution of the values $|\zeta(1+it)|$ (in [7]), and showed that uniformly in the region $\tau \leq \log_2 T-\log_3 T$, the measure of points $t\in[T,2T]$ for which $|\zeta(1+it)|\geq e^{\gamma}\tau$, equals
$$
T\exp\left(-\frac{e^{\tau-C_2}}{\tau} \left(1 +
o(1)\right)\right), \text{ where }C_2:=1+\displaystyle{\int_0^1\log I_0(t)\frac{dt}{t^2} + \int_1^{\infty} (\log
 I_0(t)-t)\frac{dt}{t^2}},\tag{2}$$
 and $I_0(t): = \sum_{n=0}^{\infty}
(t/2)^{2n}/n!^2$ is the modified Bessel function of order $0$.

Recently, Liu, Royer and Wu [19] studied the distribution of the values  of automorphic $L$-functions at $s=1$ in the weight aspect, and showed that the distribution function of $L(1,f)$ where $f$ varies over elements of $S_k^p(1)$, has the same shape as (1) and (2).

The methods used to prove these results are similar and consist mainly in two steps. The first one is to prove that large complex moments of our family of $L$-functions and those of the corresponding random model are roughly equal, and the second one is to compute these moments using techniques from analytic number theory and probability. Therefore it seems interesting to prove a very general statement for the distribution of various families of $L$-functions using these ideas. In this work, we investigate the distribution of a general class of random Euler products, satisfying very natural conditions (confer conditions 1-4 below), then we deduce information about values of ``$L$-functions'' via their corresponding random model. Among our results, we will prove a general formula for the constant involved in the corresponding distribution function.  Furthermore, we recover all the previous results proved by Granville-Soundararajan ([7] and [8]) and Liu-Royer-Wu [19]. Indeed, the corresponding random models for the families $|\zeta(1+it)|$, $L(1,\chi_d)$ and $L(1,f)$ satisfy our conditions 1-4 below.
\noindent Another novelty in our work is to realize that we only need to compute large integral moments of $L$-functions rather than the complex ones. Indeed computing large complex moments of $L$-functions becomes a hard task when a completely multiplicative structure for their coefficients lacks.

 Let $d$ be a positive integer, and $\Cal{P}$ be the set of all prime numbers. For $p\in \Cal{P}$ and $1\leq j\leq d$, let $\theta_j(p)$  be random variables distributed on $[-\pi,\pi]$ and satisfying:

\roster
 \item""{\bf Condition 1}. ${\Bbb E}(e^{i\theta_j(p)})=0$, for all $p\in \Cal{P}$ and $1\leq j\leq d$.

\item""{\bf Condition 2}. $\theta_j(p)$ and $\theta_k(q)$ are independent random variables for  $p\neq q$.

\item""{\bf Condition 3}. The random variables $X(p):=\sum_{j=1}^de^{i\theta_j(p)}/d,$ are identically distributed, for every $p\in\Cal{P}$.

\item""{\bf Condition 4}. There exists an absolute constant $\alpha>0$ such that for all primes $p$ and all $\epsilon>0$,  we have
$ \pr\left(|\theta_1(p)|\leq \epsilon,...,|\theta_d(p)|\leq \epsilon \right)\gg \epsilon^{\alpha}.$
\endroster

Consider the following random Euler products

$$ L(1,X):=\prod_{p\in \Cal{P}}\prod_{j=1}^{d}\left(1-\frac{e^{i\theta_j(p)}}{p}\right)^{-1}.$$
 Our aim  then is to study the behavior of

$$\Phi(\tau):=\pr\left(|L(1,X)|>(e^{\gamma}\tau)^d\right).$$
\proclaim {Theorem 1} Let $d$ be a positive integer. For $1\leq j\leq d$ and $p\in\Cal{P}$, let $\theta_{j}(p)$ be random variables distributed on $[-\pi,\pi]$ and satisfying conditions 1-4. For $\tau\gg 1$, we have
$$ \Phi(\tau)=\exp\left(-\frac{e^{\tau-A_X}}{\tau}\left(1+O\left(\frac
{1}{\sqrt{\tau}}\right)\right)\right),$$
where
$$ A_X:= 1+ \int_0^{\infty} \frac{f(t)}{t^2}dt, \text{ with } f(t):=\left\{\aligned &\log\ex\left(e^{\text{Re}(X)t}\right)\ \ \ \ \ \text{ if } 0\leq t<1,\\  &\log\ex\left(e^{\text{Re}(X)t}\right)-t \ \text{ if } 1\leq t, \endaligned \right.
$$
and $X$ is a random variable having the same distribution as the $X(p)$. Furthermore $A_X$ is convergent by Lemma 1.1 below.
\endproclaim

In [1], Cogdell and Michel proved that large complex moments at $s=1$ of the family of $k$-th symmetric power $L$-functions of primitive cusp forms of weight $2$ and large prime level $q$, coincide with those of the adequate probabilistic model (constructed from the Sato-Tate distribution) as $q\to\infty$, assuming the following hypothesis:

\noindent {\bf Hypothesis} Sym$^k(q)$: For all $f\in S_2^p(q)$ the $k$-th symmetric power $L$-function of $f$ is automorphic, that is it coincides with the $L$-function of a certain  cuspidal automorphic representation of $GL(k+1)/{\Bbb Q}$.

This hypothesis is predicted by the Langlands functoriality conjectures and is effectively proved for the symmetric powers up to $4$.
In view of the Petersson trace formula, it is arguably more natural to consider the weighted arithmetic distribution function
$$ \Phi_q(\text{Sym}^k,\tau)=\left(\sum_{f\in S_2^p(q)}\omega_f\right)^{-1}
\sum\Sb f\in S_2^p(q)\\ L(1,\text{Sym}^kf)\geq (e^{\gamma}\tau)^{k+1} \endSb \omega_f,$$
where $ \omega_f:=1/(4\pi\|f\|)$ is the usual harmonic weight. In their paper [1], Cogdell and Michel did not determine a good estimate for the appropriate distribution function, and noted that it will be interesting to have such an estimate. Liu, Royer and Wu [19] noted that their method should work in this case but with additional technical issues. Using Theorem 1 we get a good estimate for this distribution function. Indeed, the corresponding random model for this family satisfies conditions 1-4, since in this case we have $\theta_j(p)=(k-2j)\theta_p$ for $0\leq j \leq k$, where the $\{\theta_p\}_{p \in \Cal{P}}$ are independent random variables distributed on $[0,\pi]$ according to the Sato-Tate measure $d\mu_{st}=\frac{2}{\pi}\sin^2(\theta)d\theta$.
\proclaim {Theorem 2} Let $k\geq 1$ be an integer and $q$ be a large prime such that Hypothesis Sym$^k(q)$ holds. Then uniformly in the region $\tau\leq \log_2 q-\log_3 q-2\log_4 q$ we have
$$ \Phi_q(\text{Sym}^k,\tau)=\exp\left(-\frac{e^{\tau-A_k}}{\tau}
\left(1+O\left(\frac{1}{\sqrt{\tau}}\right)\right)\right),$$
where $A_k=1+\int_0^1\frac{h_k(t)}{t^2}dt +\int_1^{\infty}\frac{h_k(t)-t}{t^2}dt$ and
$$ h_k(t)=\log\left(\frac{2}{\pi}\int_0^{\pi}\exp\left(\frac{t}{k+1}
\sum_{j=0}^k\cos(\theta(k-2j))\right)\sin^2\theta d\theta\right).$$
\endproclaim

Another application of our work concerns the distribution of functions of the Selberg class in the $t$-aspect. This class $S$, introduced by Selberg
[26] (see also the nice survey of Kaczorowski and Perelli [12]), is the class of Dirichlet series
$$ F(s)=\sum_{n=1}^{\infty}\frac{a_F(n)}{n^s}, \text{ for Re}(s)>1,$$
satisfying the following axioms

\roster
\item""{\bf Axiom 1} Analyticity: $(s-1)^lF(s)$ is an entire function of finite order for some non-negative integer $l$.

\item""{\bf Axiom 2} Ramanujan hypothesis: $a_F(n)\ll_{\epsilon}n^{\epsilon} $ for any fixed $\epsilon>0$.

\item""{\bf Axiom 3} Functional equation: $F$ satisfies the functional equation $$ \Phi(s)=\omega\overline{\Phi(1-\overline{s})},$$
    where
$$ \Phi(s)= Q_F^s\prod_{i=1}^k\Gamma(w_is+\mu_i)F(s),$$
and $|\omega|=1$, $Q_F>0$, $w_i>0$ and Re$(\mu_i)\geq 0$ are parameters depending on $F$.

\item""{\bf Axiom 4} Euler product: $a_1=1$, and
$$ \log F(s)=\sum_{n=1}^{\infty}\frac{b_n}{n^s},$$
where $b_n=0$ unless $n=p^m$ with $m\geq 1$, and $b_n\ll n^{\lambda}$ for some $\lambda<1/2$.
\endroster

All $L$-functions appearing in number theory are conjectured to belong to the Selberg class, and most of these $L$-functions are indeed in $S$. Examples include the Riemann zeta function $\zeta(s)$, Dirichlet $L$-functions of primitive characters $L(s,\chi)$, the Dedekind zeta function $\zeta_K(s)$ of an algebraic number field $K$, the $L$-function associated with a normalized primitive cusp form of some congruence subgroup of $SL_2({\Bbb Z})$, and the $L$-function of a Rankin-Selberg convolution of two normalized primitive cusp forms. However, the $L$-functions attached to $GL(2)$ Maass cusp forms are not known to belong to $S$ since the Ramanujan hypothesis (axiom 2) is missing for these $L$-functions. Therefore, in order to include this latter case in our setting we define another class of functions which satisfy axioms 1, 3 and 4 and for which the Ramanujan  hypothesis is replaced by the following assumption:

{\bf Strong Ramanujan hypothesis on average}. There exists a constant $\beta_F\geq 0$ such that
$$ \sum_{n\leq x}|a_F(n)|\ll x(\log x)^{\beta_F}. \tag{3}$$

We denote the class of such functions by $S^{a}$. By a classical result (Theorem 8.3 of [10] combined with Cauchy-Schwarz inequality) we know that $L$-functions attached to $GL(2)$ Maass cusp forms are in $S^a$. Indeed these $L$-functions satisfy assumption (3) with $\beta_F=0$.

It is conjectured that the Selberg class consists only of automorphic $L$-functions,
and for those the Euler product (axiom 4) has the form:
$$ F(s)=\prod_{p\in \Cal{P}}\prod_{i=1}^d\left(1-\frac{\alpha_{i,F}(p)}{p^s}\right)^{-1}, \hbox{ for Re}(s)>1, \tag{4}$$
 where $\alpha_{i,F}(p)\neq 0$ for all primes except finitely many. The $\alpha_{i,F}(p)$ are complex numbers called {\it the local roots} of $F$ at $p$ and $d\in{\Bbb N}$ is called {\it the degree} of $F$. In this case we say that $F$ has a {\it polynomial Euler product}.

 Let $S^p$ be the class of functions which satisfy axioms 1-4 of $S$ and which have a polynomial Euler product. All examples of elements of $S$ cited above are in $S^p$. A classical result (see Lemma 4.1 below) shows that the Ramanujan hypothesis along with the polynomial Euler product implies that
 $$ |\alpha_{i,F}(p)|\leq 1, \text{ for all } p\in \Cal{P} \text{ and } 1\leq i\leq d,\tag{5}$$
and therefore, we can show the following key fact (see Lemma 4.1 for a proof)
$$ S^p\subset S\cap S^a,\tag{6}$$
while it is conjectured that $S^p=S.$

Let $F\in S^a$. Our aim is to study the distribution of large values of $|F(1+it)|$ for $t\in[T,2T]$, and $T$ large. To this end we construct a probabilistic model based on the totally multiplicative function $X(.)$, where the $\{X(p)\}_{p\in \Cal{P}}$ are independent random variables uniformly distributed on the unit circle $\Bbb{U}$, and for a positive integer $n=p_1^{a_1}...p_r^{a_r}$ we have $X(n):=X(p_1)^{a_1}...X(p_r)^{a_r}.$ Then define the following random series
$$ F(1,X):=\sum_{n=1}^{\infty}\frac{a_F(n)X(n)}{n}, \hbox{ (these series converge with probability }1).$$
We first prove that large integral moments of $|F(1+it)|$ and those of $|F(1,X)|$ are roughly equal, assuming that $F$ satisfies the following assumption

\noindent{\bf Weak zero density estimate near $\sigma=1$ (WZD)}. Let $$N_F(\sigma,T):=|\{ \rho=\beta+i\gamma: F(\rho)=0, \sigma\leq \beta\leq 1,  |\gamma|\leq T\}|$$
Then there exist two constants $\lambda _F>0$ and $0<\epsilon_F<1$  such that
$$ N_F(\sigma,T)\ll T^{\lambda_F(1-\sigma)}, \text{ uniformly for } 1-\epsilon_F<\sigma<1.$$
This assumption has been proved for a wide range of $L$-functions. Indeed it holds for all elements of the Selberg Class $S$ by the work of Kaczorowski and Perelli [11], and for $L$-functions attached to $GL(2)$-Maass cusp forms by the work of Sankaranarayanan and Sengupta [25]. We prove

\proclaim {Theorem 3}
Let $F\in S^a$ and satisfies (WZD). Let $T>0$ be large, and take $A>0$. Then for all positive integers $k$ in the range $1\leq k\leq \log
T/(B(\log_2T)^{\beta_F+2})$ (for a suitably large constant $B=B(A,F)$) , we have
$$\frac{1}{T}\int_T^{2T}|F(1+it)|^{2k}dt = {\Bbb E}\left(|F(1,X)|^{2k}\right)\left(1+O\left(\frac{1}{\log^A T}\right)\right).\tag{7}$$
Moreover if $F\in S^p$, then this latter asymptotic holds in the wider range $1\leq k\leq \log
T/(C\log_2T\log_3 T)$, for a suitably large constant $C=C(A,F)$.
\endproclaim
{\it Remark 1.} To improve the range of uniformity for the moments of elements of $S^p$, we use (5) to show that if $F\in S^p$ then  $\sum_{p\leq x}|a_F(p)|/p\ll \log\log x.$
Furthermore, it is conjectured that for all $F\in S$ there exists some constant $\kappa_F>0$ such that
$$ \sum_{p\leq x}\frac{|a_F(p)|}{p}=\kappa_F\log\log x+E_F(x),\tag{8}$$
where $E_F(x)$ is an error term, and we believe further that $E_F(x)=O(1)$. This conjecture is intimately related to the Deep Selberg orthogonality conjectures (see [26]). For $L$-functions of degree one ($\zeta(s)$ and $L(s,\chi)$) this is indeed true, and for $L$-functions attached to elliptic curves and their symmetric powers (up to $4$) we can show that (8) holds with $E_F(x)=o(\log\log x)$ (see remark 3 below).

Using assumption (8) and Theorem 3 we can exhibit large values of $|F(1+it)|$, when $F\in S^p$. Indeed we have
\proclaim{Corollary 1} Let $F\in S^p$ and satisfies assumption (8). Then for $T>0$ large, there exists some $t\in [T,2T]$ such that
$$|F(1+it)|\gg\left\{\aligned &(\log\log T)^{\kappa_F+o(1)}\ \ \ \ \ \text{ if } E_F(x)=o(\log\log x),\\  &(\log\log T)^{\kappa_F} \ \text{ if } E_F(x)=O(1). \endaligned \right.$$
\endproclaim
{\it Remark 2.} When $E_F(x)=O(1)$, the bound provided by Corollary 1 is best possible (up to the multiplicative constant). Indeed by a standard argument of Littlewood ([17] and [18]) one can show that, under the Generalized Riemann Hypothesis for $F$  $$|F(1+it)|\ll (\log\log t)^{\kappa_F}.$$

Therefore what remains is to study the random model closely to deduce an estimate for the distribution function of $|F(1+it)|$. First we restrict our selves to the class $S^p$. Furthermore, we believe that for most functions of $S^p$ the values $|a_F(p)|$ behave in some regular way, that is there exists some compactly supported distribution function $\psi(t)$, such that
$$ \frac{1}{\pi(x)}|\{p\leq x: |a_F(p)|\in I\}|\sim\int_{I}\psi(t)dt, \text{ as } x\to\infty, \tag{9}$$

for any interval $I$.

{\it Remark 3.} This assumption is proved to hold for the following $L$-functions:
\roster
\item"1." Degree one $L$-functions: $\zeta(s)$ and $L(s,\chi)$, in which case $\psi(t)=\delta(t-1)$, where $\delta$ is the Dirac-Delta function.
\item"2." Normalized $L$-function attached to a non-CM elliptic curve over a totally real field. In which case $\psi(t)=\frac{2}{\pi}\sqrt{1-t^2/4},$ for $t\in [0,2]$ and equals $0$ otherwise. This can be derived from the Sato-Tate measure.
\item"3." Normalized $L$-function attached to a CM elliptic curve. In this case by the work of Deuring [3] we can prove that $\psi(t)=\frac{1}{2}\delta(x)+1/(2\pi\sqrt{1-t^2/4})$, for $t\in[0,2]$ and equals $0$ otherwise.
\endroster
Moreover for our purpose we need a slightly uniform version of (9):

\proclaim {Hypothesis D}
There exists a compactly supported distribution function $\psi(t)$ (with support in some interval $[0,U]$), such that for all continuous functions $g$ we have
$$ \sum_{p\leq x}g(|a_F(p)|)=\pi(x)\left(\int_{0}^Ug(t)\psi(t)dt+o\left(\frac{1}{\log x}\right)\right), \text{ as } x\to\infty.$$
\endproclaim

Let $F\in S^{p}$ and satisfies hypothesis $D$. Let $N:=\int_0^Ut\psi(t)dt$ and $M:=\int_0^U t\log t \psi(t)dt.$ Then one can easily prove that $F$ satisfies assumption (8) with $\kappa_F=N$, and $E_F(x)=O(1).$ This implies that the following Euler product is convergent
$$ b_F:= \prod_{p\in \Cal{P}}\max_{t\in[-\pi,\pi]}\left|\prod_{i=1}^d\left(1-\frac{e^{it}\alpha_{i,F}(p)}
{p}\right)^{-1}\right|\left(1-\frac{1}{p}\right)^{N}.$$
Now define
$$ \Phi_F(\tau):=\frac{1}{T}\text{meas}\{ t\in [T,2T]: \
\ |F(1+it)| >b_F\left(e^{\gamma}\tau\right)^{N}\}.\tag{10}$$
\proclaim {Theorem 4} Let $T>0$ be large. Let $F\in S^{p}$, and satisfies Hypothesis D (with distribution function $\psi$). Then uniformly in the region $\tau\leq \log_2T-\log_3 T-2\log_4 T $, we have
$$\Phi_F(\tau)=\exp\left(-\frac{e^{\tau-C_2-M/N+\log N}}{\tau}\left(1+
o(1)\right)\right).$$

\endproclaim

In equations (1), (2) and in Theorem 2, the constant involved in the distribution function depends only on the corresponding random model. However in the case of  Theorem 4, this constant contains two parts, the first one $C_2$ depends on the random model, and the second one $M/N-\log N$ depends on the function $\psi$ which governs the distribution of the $|a_F(p)|$. This latter constant equals $0$ for $L$-functions of degree one, but is non-vanishing for $L$-functions attached to elliptic curves, and can be computed using remark 3 above.  Therefore it seems that assuming Hypothesis D (or at least assumption (9)) is necessary if one wants to know the distribution of large values of $F(1+it)$.

Let $F\in S$, and $\chi$ a primitive character mod $q$. The twist of $F$ by $\chi$ is defined for Re$(s)>1$ by the Dirichlet series
$$ F\otimes \chi(s)=\sum_{n=1}^{\infty}\frac{a_{F}(n)\chi(n)}{n^s}.$$
For general $F$ we do not know if $F\otimes\chi$ belongs to $S$. Therefore we restrict our selves to the automorphic world where this holds.
Let $\pi$ be an auto-dual cuspidal automorphic representation of $GL(m)/{\Bbb Q}$. The normalized  $L$-function of $\pi$ (so that the critical strip is $0\leq \text{Re}(s)\leq 1$) satisfies axioms 1, 3 and 4, and has a polynomial Euler product as in (4) with degree $m$. Let $L(\pi\otimes\chi_d,s)$ be the twist of $L(\pi,s)$ by the unique real quadratic character mod $d$, where $d$ is a fundamental discriminant.
The last part of our work concerns  the study of the distribution of large values of $L(\pi\otimes \chi_d,1)$ for fundamental discriminants $|d|\leq x$, where $x$ is large. We prove the analogue of Theorem 3 in this case (see Theorem 5.2) assuming that $\pi$ satisfies assumption (3) and that $|a_{\pi}(p)|\ll p^{\theta} \text{ for some } 0<\theta<1/4.$ In particular this holds for $GL(2)$ Maass cusp forms, where one can take $\theta=1/9$ by the work of Kim and Shahidi [13]. Moreover if $L(\pi,s)$ satisfies the Ramanujan hypothesis and hypothesis $D$ (with distribution function $\psi$), then we have the analogue of Theorem 4 in this case. Indeed, similarly to (10) let $\Phi_{\pi}(\tau)$ denotes the proportion of fundamental discriminants $d$ with $|d|\leq x$, such that $L(\pi\otimes \chi_d,1)>b_{\pi}(e^{\gamma}\tau)^N,$ where $ N=\int_0^U t\psi(t) dt$, $M=\int_0^U t\log t \psi(t)dt,$ and $$b_{\pi}:= \prod_{p\in \Cal{P}}\max_{\delta\in\{-1,1\}}\prod_{i=1}^d\left(1-\delta\frac{\alpha_{i,F}(p)}
{p}\right)^{-1}\left(1-\frac{1}{p}\right)^N.$$
 \proclaim{ Theorem 5}  Let $\pi$ be an auto-dual cuspidal automorphic representation of $GL(m)/{\Bbb{Q}}$ satisfying the Ramanujan hypothesis, and such that the $a_{\pi}(p)$ satisfy Hypothesis D (with distribution function $\psi$). Let $x>0$ be large. Then uniformly in the region $\tau\leq \log_2x-\log_3 x-2\log_4 x $, we have
$$\Phi_{\pi}(\tau)=\exp\left(-\frac{e^{\tau-C_1-M/N+\log N}}{\tau}\left(1+
o(1)\right)\right).$$

\endproclaim

\head {1. Distribution of Random Euler Products: Proof of Theorem 1}\endhead

For $1\leq j\leq d$ and $p\in\Cal{P}$, let $\theta_{j}(p)$ be random variables distributed on $[-\pi,\pi]$ and satisfying conditions 1-4. Let $X$ be a random variable having the same distribution as the $X(p)$, define $h(t):=\log\ex\left(e^{\text{Re}(X)t}\right)$, and let $f(t)$ be as in Theorem 1. First we prove some useful properties of the function $f$.
\proclaim{ Lemma 1.1}

 i) $f$ satisfies the growth conditions $$
  \ \ \ f(t)=\left\{\aligned &O\left(t^2\right) \ \ \ \ \text{ if } 0\leq t<1\\  &O\left(\log t+1\right) \ \text{ if } 1\leq t. \endaligned \right.
$$

ii) $h$ and $f$ are  continuous differentiable Lipshitz functions on $[0,+\infty)$.

iii) $A_X$ is convergent.
\endproclaim

\demo{ Proof}

i) First let $0\leq t<1$. Noting that $\ex(\text{Re}(X))=0$ by condition 1, and using the Taylor expansion of $e^{\text{Re}(X)t}$, we deduce that $f(t)=\log(1+O(t^2))=O(t^2).$
Now for $t\geq 1$ , we have that $f(t)= \log\ex\left(e^{(\text{Re}(X)-1)t}\right).$ Since $\text{Re}(X)\leq 1$, then $f(t)\leq 0$. Let $p$ be  prime. By condition 3 we have
$$ e^{f(t)}\geq \pr\left(|\theta_1(p)|\leq \epsilon,...,|\theta_d(p)|\leq \epsilon \right)e^{(\cos\epsilon-1)t}, \text{ for all } \epsilon>0.$$
Choose $\epsilon=1/t$, and use condition 4 to get the result.

ii) By  the Taylor expansion of   $e^{\text{Re}(X)t}$ and the fact that $|\text{Re} (X)|\leq 1$,  $h$ is smooth and
$$ |h'(t)|=\frac{\left|\ex\left(\text{Re}(X)e^{\text{Re}(X)t}\right)\right|}{\ex\left(e^{\text{Re}(X)t}\right)}\leq \frac{\ex\left(\left|\text{Re}(X)\right|e^{\text{Re}(X)t}\right)}{\ex\left(e^{\text{Re}(X)t}\right)}\leq 1.$$
Hence $h$ is Lipschitz with constant $1$. This implies also that $f$ is differentiable and Lipshitz with constant $2$.

iii) This follows from (i).
\enddemo

Now we are ready to compute the real moments of the random variable $|L(1,X)|$. We prove the following proposition

\proclaim{ Proposition 1.2} Let $r>d^{8}$. Then we have
$$ \log\ex\left(|L(1,X)|^r\right)=rd\log_2(rd)+rd\gamma+ \frac{rd}{\log rd}\left(A_X-1+O\left(\frac{1}{\log r}\right)\right).$$
\endproclaim

\demo{ Proof } We have that
$$ |L(1,X)|^r=\prod_{p\in\Cal{P}}\prod_{j=1}^d
\left|1-\frac{e^{i\theta_j(p)}}{p}\right|^{-r}= \prod_{p\in\Cal{P}}\prod_{j=1}^d
\left(1-\frac{2\cos\theta_j(p)}{p}+\frac
{1}{p^2}\right)^{-r/2}.
$$
Let $$E_p:=\ex\left(\prod_{j=1}^d
\left(1-\frac{2\cos\theta_j(p)}{p}+\frac
{1}{p^2}\right)^{-r/2}\right).$$ Then by condition 2 we know that $\ex(|L(1,X)|^r)=\prod_{p\in \Cal{P}}E_p$.

\noindent {\bf Case 1}.  $p>\sqrt{rd}$. In this case
$$
\prod_{j=1}^d
\left(1-\frac{2\cos\theta_j(p)}{p}+\frac
{1}{p^2}\right)^{-r/2}
=\exp\left(\frac{rd}{p}\text{Re}(X(p))+
O\left(\frac{rd}{p^2}\right)\right).\tag{1.1}
$$
{\bf Case 2}. $p\leq \sqrt{rd}$. Let $$E_p*:=E_p\left(1-\frac{1}{p}\right)^{rd}.$$ Then we have
$$ E_p*=\ex\left(\prod_{j=1}^d\left(1-\frac{2(\cos\theta_j(p)-1)}
{p(1-1/p)^2}\right)^{-r/2}\right)\leq 1.\tag{1.2}$$
Moreover for all $\epsilon>0$ we have that
$$ E_p*\geq \pr\left(|\theta_1(p)|\leq \epsilon,...,|\theta_d(p)|\leq \epsilon \right)\prod_{j=1}^d\left(1-\frac{2(\cos\epsilon-1)}
{p(1-1/p)^2}\right)^{-r/2}.$$
Choosing $\epsilon=1/(dr)$ and using condition 4, we get
$$ E_p*\gg \frac{1}{(dr)^{\alpha}}\exp\left(-\frac{1}{2rd^2p(1-1/p)^2}
\right)\gg\frac{1}{(dr)^{\alpha}}.\tag{1.3}$$
 Therefore by (1.2) and (1.3) we get
$$ \log E_p=-rd\log \left(1-\frac{1}{p}\right)+O(\log r).$$
Hence by (1.1) we conclude that
$$
\sum_{p\in \Cal{P}}\log E_p=-rd\sum_{p\leq \sqrt{rd}}\log \left(1-\frac{1}{p}\right)
+ \sum_{\sqrt{rd}<p}h\left(\frac{rd}{p}\right)+O\left(\sum_{p\leq \sqrt{rd}}\log r+\sum_{\sqrt{rd}<p}\frac{dr}{p^{2}}\right).
$$
For the error terms, since $r>d^8$ we have
$$ \sum_{\sqrt{rd}<p}\frac{dr}{p^{2}}+\sum_{p\leq \sqrt{rd}}\log r\ll r^{3/4}.$$
Therefore using that
$$ \sum_{\sqrt{rd}\leq p\leq rd}\log\left(1-\frac 1p\right)+\frac 1p\ll (rd)^{-1/2},$$ we get
$$ \sum_{p\in \Cal{P}}\log E_p= -rd\sum_{p\leq rd}\log\left(1-\frac{1}{p}\right)+\sum_{\sqrt{rd}<p}f
\left(\frac{rd}{p}\right)+O\left(r^{3/4}\right).\tag{1.4}$$
Let $r_1=r^{3/2}.$ Then by Lemma 1.1
$$ \sum_{p>r_1}f\left(\frac{rd}{p}\right)\ll \sum_{p>r_1}\frac{(rd)^2}{p^2}\ll r^{3/4}.\tag{1.5}$$
To evaluate the sum over $f(rd/p)$ we use the prime number theorem
 in the form
$$\pi(t)=\int_2^t\frac{du}{\log u}+O\left(te^{-8\sqrt{\log t}}\right).$$
Therefore
$$ \sum_{\sqrt{rd}<p\leq r_1}f\left(\frac{rd}{p}\right)=
\int_{\sqrt{rd}}^{r_1}f\left(\frac{rd}{t}\right)d\pi(t)=
\int_{\sqrt{rd}}^{r_1}f\left(\frac{rd}{t}\right)\frac{dt}{\log t}+E_1,\tag{1.6}$$
where
$$
E_1\ll f\left(\sqrt{rd}\right)\sqrt{rd}e^{-4\sqrt{\log r}}+f\left(\frac{rd}{r_1}\right)r_1e^{-4\sqrt{\log r}}+\int_{\sqrt{rd}}^{r_1}\frac{rd}{t^2}\left|f'\left(\frac{rd}{t}\right)\right|te^{-8\sqrt{\log t}}dt.$$
Now by Lemma 1.1 we can see that $E_1\ll rde^{-4\sqrt{\log r}}.$
To estimate the main term we make the change of variables $T=rd/t$. Hence we have
$$
\int_{\sqrt{rd}}^{r_1}f\left(\frac{rd}{t}\right)\frac{dt}{\log t}=
rd\int_{dr^{-1/2}}^{\sqrt{rd}} \frac{f(T)}{T^2\log(rd/T)}dT.\tag{1.7}
$$
In the range $dr^{-1/2}\leq t\leq \sqrt{rd}$, we have
$$ \frac{1}{\log(rd/t)}=\frac{1}{\log(rd)}\frac{1}{1-\frac{\log t}{\log(rd)}}=\frac{1}{\log(rd)}+O\left(\frac{\log t}{\log^2 r}\right).$$
Therefore
$$ \int_{dr^{-1/2}}^{\sqrt{rd}} \frac{f(t)}{t^2\log(rd/t)}dt=\frac{1}{\log (rd)} \int_{dr^{-1/2}}^{\sqrt{rd}} \frac{f(t)}{t^2}dt+O\left(\frac{1}{\log^2r}\right),\tag{1.8}$$
using that
$$ \int_{0}^{
\infty}\frac{f(t)\log(t)}{t^2}dt\ll 1,$$
which follows from Lemma 1.1.
Finally by the same Lemma we deduce that
$$ \int_{dr^{-1/2}}^{\sqrt{rd}} \frac{f(t)}{t^2}dt=A_X-1+O\left(\frac{d}{\sqrt{r}}+\frac{\log rd}{\sqrt{rd}}\right).\tag{1.9}$$
Therefore from  equations (1.4)-(1.9) we deduce that
$$ \sum_{p\in\Cal{P}}\log E_p= -rd\sum_{p\leq rd}\log\left(1-\frac{1}{p}\right)+\frac{rd}{\log rd}(A_X-1)+O\left(\frac{rd}{\log^2 r}\right).$$
Finally by the prime number theorem we have
$$-\sum_{p\leq rd}\log\left(1-\frac{1}{p}\right)=\log_2(rd)+\gamma+O\left(e^{-2\sqrt{\log r}}\right),\tag{1.10}$$ which implies the result.
\enddemo

\demo{ Proof of Theorem 1} Let the $\{X(p)\}_{p\in \Cal{P}}$ be defined on a probability space $(\Omega,\mu).$
Then
$$
\align
&rd\int_0^{\infty}\Phi(t)t^{rd-1}dt =
rd\int_0^{\infty}t^{rd-1}\int_{|L(1,X(\omega))|>(e^{\gamma}t)^d}d\mu(\omega)dt\\
&=\int_{\Omega} \frac{|L(1,X(\omega))|^r}{(e^{\gamma rd })}d\mu(\omega)=\ex\left(|L(1,X)|^r\right)e^{-\gamma rd}.\\
\endalign
$$
Therefore by Proposition 1.2 we deduce that
$$ \int_0^{\infty}\Phi(t)t^{rd-1}dt=(\log rd)^{rd}\exp\left(\frac{rd}{\log rd}\left(A_X-1+O\left(\frac{1}{\log r}\right)\right)\right).\tag{1.11}$$
Let $\tau=\log rd +A_X$, and $R=r e^{\delta}$ where $\delta>0$ will be chosen later. By (1.11)  we have
$$
\align
\int_{\tau+\delta}^{\infty}\Phi(t)t^{rd-1}dt &\leq (\tau+\delta)^{rd\left(1-e^{\delta}\right)}\int_0^{\infty}\Phi(t)t^{Rd-1}dt\\
&=(\tau+\delta)^{rd\left(1-e^{\delta}\right)}(\log rd+\delta)^{rde^{\delta}}\exp\left(\frac{rde^{\delta}}{\log rd}\left(A_X-1+O\left(\frac{1}{\log r}\right)\right)\right)\\
&=(\log rd)^{rd}\exp\left(\frac{rd}{\log rd}\left(A_X-1+(\delta+1-e^{\delta})+O\left(\frac{1}{\log r}\right)\right)\right).\tag{1.12}\\
\endalign
$$
Now choose $\delta=c/\sqrt{\log r},$ where $c$ is a suitably large constant. Then by (1.11) and (1.12) we deduce that
$$
\int_{\tau+\delta}^{\infty}\Phi(t)t^{rd-1}dt\leq \left(\int_0^{\infty}\Phi(t)t^{rd-1}dt\right) \exp\left(-\frac{r}{\log^2 r}\right).\tag{1.13}$$
A similar argument shows that
$$
\int_{0}^{\tau-\delta}\Phi(t)t^{rd-1}dt\leq \left(\int_0^{\infty}\Phi(t)t^{rd-1}dt\right) \exp\left(-\frac{r}{\log^2 r}\right).\tag{1.14}$$
Therefore by (1.11), (1.13) and (1.14) we deduce that
$$ \int_{\tau-\delta}^{\tau+\delta}\Phi(t)t^{rd-1}dt=(\log rd)^{rd}\exp\left(\frac{rd}{\log rd}\left(A_X-1+O\left(\frac{1}{\log r}\right)\right)\right).\tag{1.15}$$
Since $\Phi(t)$ is a non-increasing function we have that
$$ \Phi(\tau+\delta)\tau^{rd}\exp\left(O\left(\frac{\delta r}{\tau}\right)\right)\leq \int_{\tau-\delta}^{\tau+\delta}\Phi(t)t^{rd-1}dt\leq \Phi(\tau-\delta)\tau^{rd}\exp\left(O\left(\frac{\delta r}{\tau}\right)\right).$$
This along with (1.15) imply that
$$ \Phi(\tau+\delta)\leq \exp\left(-\frac{e^{\tau-A_X}}{\tau}\left(1+O\left(\delta\right)\right)\right)\leq \Phi(\tau-\delta),$$
and thus
$$ \Phi(\tau)=\exp\left(-\frac{e^{\tau-A_X}}{\tau}\left(1+O\left(\frac{1}{\sqrt{\tau}}\right)\right)\right),$$
completing the proof.
\enddemo

\head {2. Distribution of Symmetric power $L$-functions: Proof of Theorem 2}\endhead

In this section we study the distribution of the family $L(1,\text{Sym}^kf)$, for $f\in S_2^p(q)$, where $q$ is a large prime number. Given a sequence $(\alpha_f)_{f\in S_2^p(q)}$, its harmonic average is defined as the sum
$$ \sumh_{f\in S_2^p(q)}\alpha_f=\sum_{f\in S_2^p(q)}\frac{\alpha_f}{4\pi\|f\|},$$
and if $S\subset S_2^p(q)$ then we will let $|S|_h$ denote the harmonic measure of $S$, that is
$$ |S|_h:=\sumh_{f\in S}1.$$
Such averaging is natural in view of the two following facts
$$ |S_2^p(q)|_h=1+O\left(\frac{\log q}{q^{3/2}}\right), \text{ and } \frac{1}{q(\log q)^3}\ll \omega_f\ll \frac{\log q}{q};\tag{2.1}$$
so that the harmonic weight $\omega_f$ is not far from the natural weight $1/|S_2^p(q)|$ (since $|S_2^p(q)|\asymp q$), and it defines asymptotically a probability measure on $S_2^p(q)$.

To describe the corresponding random model for the family $L(1,\text{Sym}^kf)$, we consider the compact group $G=SU(2)$ endowed with its natural Haar measure $\mu_G$; we then let $G^{\natural}$ be the set of conjugacy classes of $G$ endowed with the Sato-Tate measure $\mu_{st}$ (i.e. the direct image of $\mu_G$ by the canonical projection). By Weyl's integration formula, the map
$$ \theta\rightarrow g(\theta)^{\natural}= {\pmatrix e^{i\theta}&0\\0&e^{-i\theta}\endpmatrix}^{\natural},$$
identifies $G^{\natural}$ with the interval $[0,\pi]$ and $\mu_{st}$ with the distribution
$$d\mu_{st}(\theta)=\frac{2}{\pi}\sin^2(\theta)d\theta.$$
Consider a probability space $(\Omega,\mu)$, and let $\{g^{\natural}_p\}_{p\in \Cal{P}}$ be a sequence of independent random variables, with values in $G^{\natural}$ and distributed according to the measure $\mu_{st}$. Consider the following random Euler product
$$ L(1,\text{Sym}^kg^{\natural}):= \prod_{p\in \Cal{P}}\det\left(I-p^{-1}\text{Sym}^kg^{\natural}_p\right)^{-1}, $$
which converges with probability $1$. Cogdell and Michel [1] proved that large complex harmonic moments of $L(1,\text{Sym}^kf)$ and those of the random product $L(1,\text{Sym}^kg^{\natural})$ are roughly equal. More precisely they showed that

\proclaim{ Theorem 2.1 (Theorem 1.3 of [1])} Let $k\geq 1$ be an integer and $q$ a prime such that Hypothesis $\text{Sym}^k(q)$ holds. Then there exists $C=C(k)>0$ and $\delta=\delta(k)>0$ and an exceptional set $S_{2,ex}^p(q)\subset S_2^p(q)$ with at most one element such that, for any complex number $z$ satisfying $|z|\leq C\log q/(\log_2 q\log_3 q)$, we have
$$ \frac{1}{|S_2^p(q)\backslash S_{2,ex}^p(q)|_h}\sumh_{f\in S_2^p(q)\backslash S_{2,ex}^p(q)}L(1,\text{Sym}^kf)^{z}=
\ex\left(L(1,\text{Sym}^kg^{\natural})^z\right)+ O_k\left(\exp\left(-\delta\frac{\log q}{\log_2 q}\right)\right),$$
the implied constant depending on $k$ only.

\endproclaim

Let $$\Phi_q^{\star}(\text{Sym}^k,\tau)=\frac{1}{|S_2^p(q)\backslash S_{2,ex}^p(q)|_h}
\sumh\Sb f\in S_2^p(q)\backslash S_{2,ex}^p(q)\\ L(1,\text{Sym}^kf)\geq (e^{\gamma}\tau)^{k+1} \endSb 1.$$
Removing at most one exceptional element from $S_2^p(q)$ does not affect the distribution function. Indeed by (2.1) we can see that
$$  \Phi_q(\text{Sym}^k,\tau)=\Phi_q^{\star}(\text{Sym}^k,\tau)+O\left(\frac{\log q}{q}\right).\tag{2.2}$$
Now we have to compute the moments of the random Euler product $L(1,\text{Sym}^kg^{\natural})$. Let $r>0$, by Weyl's integration formula we have
$$
\align
\ex\left(L(1,\text{Sym}^kg^{\natural})^r\right)&=\prod_{p\in \Cal{P}}\int_{G^{\natural}}\det
\left(I-p^{-1}\text{Sym}^kg^{\natural}\right)^{-r}d\mu_{st}(g^{\natural})\\
&=\prod_{p\in \Cal{P}}\left(\frac{2}{\pi}\int_{0}^{\pi}\det
\left(I-p^{-1}\text{Sym}^kg(\theta)\right)^{-r}\sin^2\theta d\theta\right).\\
\endalign
$$
Since $\text{Sym}^kg(\theta)$ is a diagonal matrix with diagonal entries $(e^{i\theta(k-2j)})_{0\leq j\leq k}$, then
$$\det
\left(I-p^{-1}\text{Sym}^kg(\theta)\right)= \prod_{j=0}^k\left(1-\frac{e^{i\theta(k-2j)}}{p}\right).$$

 Let $\{\theta_p\}_{p\in\Cal P}$ be independent random variables taking values on $[0,\pi]$ and endowed with the Sato-Tate measure $(2/\pi)\sin^2\theta d\theta$, and let $\theta_j(p)=(k-2j)\theta_p$ for $0\leq j\leq k$. Define the random variables $X_{st}(p):=\sum_{j=0}^k\theta_j(p)/(k+1),$  and the random Euler product $L(1,X_{st})$
 $$ L(1,X_{st}):=\prod_{p\in \Cal{P}}\prod_{j=0}^k\left(1-\frac{e^{i\theta_p(k-2j)}}{p}\right)^{-1} .$$ The local factors of this random Euler product equal
 $$ \left(1-\frac{1}{p}\right)^{-\delta_{2|k}}\prod_{0\leq j<k/2}\left(1-\frac{2\cos(\theta_p(k-2j))}{p}+\frac{1}{p^2}\right)^{-1},$$ where $\delta_{2|k}=1$ if $k$ is even and $0$ otherwise. This implies that $|L(1,X_{st})|=L(1,X_{st}).$ Also one can easily see that the $\theta_j(p)$ satisfy conditions 1-4 of Theorem 1 and Proposition 1.2. and that
$$ \ex\left(L(1,\text{Sym}^kg^{\natural})^r\right)= \ex\left(L(1,X_{st})^r\right).\tag{2.3}$$
Now we are ready to prove Theorem 2

\demo{Proof of Theorem 2} Let $r>0$ be large. Then  by Theorem 2.1 and (2.3) we have that uniformly for
$r\leq C\log q/(\log_2 q\log_3 q)$

$$
\align
&r(k+1)\int_0^{\infty}\Phi_q^{\star}(\text{Sym}^k,t)t^{r(k+1)-1}dt=\frac{r(k+1)}{|S_2^p(q)\backslash S_{2,ex}^p(q)|_h}\int_0^{\infty}t^{r(k+1)-1}
\sumh\Sb f\in S_2^p(q)\backslash S_{2,ex}^p(q)\\ L(1,\text{Sym}^kf)\geq (e^{\gamma}t)^{k+1} \endSb 1 dt\\
&=\frac{e^{-\gamma(k+1)r}}{|S_2^p(q)\backslash S_{2,ex}^p(q)|_h} \sumh_{f\in S_2^p(q)\backslash S_{2,ex}^p(q)}L(1,\text{Sym}^kf)^r=e^{-\gamma(k+1)r}\ex\left(L(1,X_{st})^r\right)+
o\left(1\right).\\
\endalign
$$
Since $L(1,X_{st})$ is real and positive, we deduce by Proposition 1.2 that
$$ \ex\left(L(1,X_{st})^r\right)=e^{\gamma(k+1)r}(\log r(k+1))^{r(k+1)}
\exp\left( \frac{r(k+1)}{\log (r(k+1))}(A_k-1)+O\left(\frac{r}{\log^2r}\right)\right), $$
where $A_k$ is defined as in Theorem 2. Therefore
$$ \int_0^{\infty}\Phi_q^{\star}(\text{Sym}^k,t)t^{r(k+1)-1}dt=(\log r(k+1))^{r(k+1)}
\exp\left( \frac{r(k+1)}{\log (r(k+1))}\left(A_k-1+O\left(\frac{1}{\log r}\right)\right)\right),$$
and thus by the same saddle point method used to prove Theorem 1 (taking $\tau=\log(r(k+1))+A_k$), we deduce that uniformly for  $\tau\leq \log_2 q-\log_3 q-2\log_4 q$, we have $$\Phi_q^{\star}(\text{Sym}^k,\tau)=\exp\left(-\frac{e^{\tau-A_k}}{\tau}
\left(1+O\left(\frac{1}{\sqrt{\tau}}\right)\right)\right).$$
Hence by (2.2) the result follows.
\enddemo

\head {3. Another class of random models } \endhead

\noindent For some families of $L$-functions, the corresponding probabilistic model has the shape
$$L(1,X)=\prod_{p\in \Cal{P}}\prod_{j=1}^{d}\left(1-\frac{X_j(p)a_j(p)}{p^s}\right)^{-1}$$ where the $X_j(p)$ are random variables distributed on the unit circle $\Bbb{U}$, and $a_j(n)$ are completely multiplicative functions. In general it is difficult to study the distribution of $|L(1,X)|$ if any of the functions $a_j(p)$ is negative or have complex values, since in this case it is difficult to find the maximum of the local factors of $|L(1,X)|$, if the $X_j(p)$ are not independent for the same prime $p$.  However we still get results in some cases where
$$ X_1(p)=X_2(p)=...=X_d(p)=X(p), \text{ for all primes } p.$$
In the first case we take the $\{X(p)\}_{p\in \Cal P}$ to be independent random variables uniformly distributed on $\Bbb{U}$ (this will correspond to the random models constructed to study functions of the Selberg class in the $t$-aspect). While in the second case the $\{X(p)\}_{p\in \Cal P}$ will be random variables taking the values $-1$ and $1$ with equal probabilities $1/2$ (this will be useful in the study of quadratic twists of $GL(m)$-automorphic $L$-functions at $s=1$). Since the analysis for the two cases is quite similar, we only describe the first one.

Let $a_1(n),...,a_d(n)$ be completely multiplicative functions such that for primes $p$ we have $ |a_i(p)|\leq 1,$
and define $a_n$  by
$$ \prod_{p\in \Cal{P}}\prod_{j=1}^{d}\left(1-\frac{a_j(p)}{p^s}\right)^{-1}=
\sum_{n=1}^{\infty}\frac{a_n}{n^s}, \text{ for Re}(s)>1, $$
so that $a_n$ is  multiplicative  and $|a_p|=|\sum_{j=1}^da_j(p)|\leq d$.
Let $\{X_u(p)\}_{p\in \Cal{P}}$ be independent random variables uniformly distributed on $\Bbb{U}$. In this section we study the distribution of the following random Euler product
$$ L(1,X_u):=\prod_{p\in \Cal{P}}\prod_{j=1}^{d}\left(1-\frac{X_u(p)a_j(p)}{p}\right)^{-1}.\tag{3.1}$$
To this end we have to estimate its real moments. Define
$$ L_p(t):=\prod_{j=1}^{d}\left(1-\frac{e^{it}a_j(p)}{p}\right)^{-1}, \text{ for } t\in[-\pi,\pi].$$
Since $|L_p(t)|$ is a continuous function of $t$, it attain  its maximum at some $\phi_p\in [-\pi,\pi].$ Let $h(t)=\log I_0(t)$. We first prove

\proclaim { Lemma 3.1} Let $r>0$ be large, and put $r_1:=r/(\log^2 r)$, and $r_2:=r\log_2 r$. Then we have the following estimates for $\log\ex\left(|L(1,X_u)|^r\right)$
$$ \log\ex\left(|L(1,X_u)|^r\right) =\left\{\aligned &r\sum_{p\leq  r_1}\log |L_p(\phi_p)| +    \sum_{r_1<p\leq r_2}h\left(\frac{r|a_p|}{p}\right)+o\left(\frac{r}{\log r}\right),\\
& r\sum_{p\leq  r}\frac{|a_p|}{p}+O(r).\\
\endaligned \right. $$
\endproclaim

\demo{Proof }
The proof goes through the lines of the proof of Proposition 1.2.
Write $\ex(|L(1,X_u)|^r)=\prod_{p\in\Cal{P}}E_p$, where
$$ E_p= \ex\left(\prod_{j=1}^{d}\left|1-\frac{X_u(p)a_j(p)}{p}\right|^{-r}\right)
=\ex\left(\prod_{j=1}^{d}\left(1-\frac{2\text{Re}(X_u(p)a_j(p))}{p}+\frac{|a_j(p)|^2}{p^2}\right)^{-r/2}\right).$$
Let $\beta_p=\arg a_p.$ In the range $p>r_1$, we have
$$
\align E_p&=\ex\left(\exp\left(\frac{r}{p}\text{Re}(a_pX_u(p))\right)\right)
\left(1+O\left(\frac{r}{p^{2}}\right)\right)\\
&=\frac{1}{2\pi}\int_{-\pi}^{\pi}
\exp\left(\frac{r|a_p|}{p}\cos(t+\beta_p)\right)dt \left(1+O\left(\frac{r}{p^{2}}\right)\right)\\
&=\exp \left(h\left(\frac{r|a_p|}{p}\right)\right)
\left(1+O\left(\frac{r}{p^{2}}\right)\right).\tag{3.2}\\
\endalign$$
Now if $p\leq r_1$, we let
$$ E_p^{\star}:=E_p|L_p(\phi_p)|^{-r}\leq 1.$$
Further $E_p^{\star}$ equals
$$ \ex\left(\prod_{j=1}^d\left(1-\frac{2\text{Re}(a_j(p)(X_u(p)-e^{i\phi_p}))}
{p\left(1-2\text{Re}(e^{i\phi_p}a_j(p))/p+|a_j(p)|^2/p^2\right)} \right)\right)^{-r/2}.$$
Let $\epsilon=\frac{1}{r^2}$. Then
$$
E_p^{\star} \geq \frac{1}{2\pi}\int_{\phi_p-\epsilon}^{\phi_p+\epsilon} \prod_{j=1}^d\left(1-\frac{2\text{Re}(a_j(p)(e^{it}-e^{i\phi_p}))}
{p\left(1-2\text{Re}(e^{i\phi_p}a_j(p))/p+|a_j(p)|^2/p^2\right)} \right)^{-r/2}dt,
$$
from which we can deduce that $E_p^{\star}\geq 1/(2\pi r^2).$ Which implies that
$$ \log E_p=r\log |L_p(\phi_p)|+O\left(\log r\right).$$
Therefore we deduce from these two cases that
$$
\sum_{p\in\Cal{P}}\log E_p=r\sum_{p\leq  r_1}\log |L_p(\phi_p)|
+\sum_{r_1<p}h\left(\frac{r|a_p|}{p}\right)+O\left(r\sum_{p>r_1}\frac{1}{p^2}+\pi(r_1)\log r\right)
.$$
The error term above is $\ll r/(\log^2 r).$
Moreover by the fact that $h(t)\ll t^2$ for $0\leq t\leq 1$, we have
$$ \sum_{p>r_2}h\left(\frac{r|a_p|}{p}\right)\ll \sum_{p>r_2}\frac{r^2|a_p|^2}{p^2}\ll \frac{r^2}{r_2\log r_2}=o\left(\frac{r}{\log r}\right),$$
which implies the first estimate. For the second one we first prove that
$$\sum_{r_1<p\leq r_2}h\left(\frac{r|a_p|}{p}\right)=o(r).$$
To this end we divide this sum into two parts, the first one $R_1$ will be over those primes for which $|a_p|\leq \delta$, and the second one $R_2$ will be over the remaining terms, where $\delta$ is a small positive number to be chosen later. Since $h(t)\ll \max(t,t^2)$ we have
$$ R_1\ll \sum_{r_1<p<r_2}\max\left\{\frac{r\delta}{p}, \frac{r^2\delta^2}{p^2}\right\}=\sum_{r_1<p<\delta r}\frac{r^2\delta^2}{p^2}+\sum_{r\delta<p<r_2}\frac{r\delta}{p}\ll \delta^2r\log r+r\frac{\log\delta}{\log r}+r\frac{\log_2r}{\log r}.$$
Moreover we know that $h(t)\ll t$ for $t\geq 1$, and that $|a_p|\leq d$. Then
$$ R_2\ll \sum_{r_1<p<r\delta}\frac{r|a_p|}{p}+\sum_{r\delta<p<r_2}
\max\left\{\frac{rd}{p}, \frac{r^2d^2}{p^2}\right\}\ll r\frac{\log\delta}{\log r}+r\frac{\log_2r}{\log r}+ \frac{r}{\delta\log(r\delta)}.$$
Now we choose $\delta=(\log r)^{-2/3}$, to deduce that
$$ R_1+R_2\ll \frac{r}{(\log r)^{1/3}}.$$ Finally our estimate follows from the two following facts
$$|L(\phi_p)|=\exp\left(\frac{|a_p|}{p}+O\left(\frac{1}{p^{2}}\right)\right), \text{ and } \sum_{r_1<p<r}\frac{|a_p|}{p}=o(1).$$
\enddemo
In order to evaluate the main term of the first estimate in Lemma 3.1, we assume further that the $a_p$ satisfy Hypothesis D with distribution function $\psi$. Let $N$ and $M$ be as in Theorem 4. To estimate the first sum we use the following convergent product
  $$b_u:=\prod_{p\in \Cal{P}}|L_p(\phi_p)|\left(1-\frac{1}{p}\right)^N.$$
  Indeed this product converges upon using Hypothesis D to show that $$\sum_{p\leq y}\frac{|a_p|}{p}= N\sum_{p\leq y}\frac 1p +O(1).$$
We now prove

 \proclaim{ Proposition 3.2} Assume that the $\{a_p\}_{p\in \Cal P}$ satisfy Hypothesis D. Let $M$ and $N$ be as in Theorem 4. For $r>0$ large enough we have
$$ \log\ex\left(|L(1,X_u)|^r\right)=rN\log_2(r)+r(\gamma N+\log b_u)+ \frac{r}{\log r}(M+N(C_2-1))+o\left(\frac{r}{\log r}\right).$$

\endproclaim

\demo{Proof}
Let $r_1$ and $r_2$ be as in Lemma 3.1. The first step is to estimate $\sum_{r_1<p<r_2}h(r|a_p|/p).$ To do so divide the interval $(r_1,r_2)$ into small interval $(u_i,u_{i+1})$ for $1\leq i\leq I$, where
$$u_1=r_1, \text{ and } u_{i+1}=u_i\left(1+\frac{1}{\log_2^2r}\right),$$
and $I$ is the smallest integer such that $ u_I>r_2$. A simple estimate shows that $I\ll (\log_2 r)^{3}.$
Now for  $1\leq i\leq I-1$ we have
$$ \sum_{u_i\leq p\leq u_{i+1}}h\left(\frac{r|a_p|}{p}\right)=\sum_{u_i\leq p\leq u_{i+1}}h\left(\frac{r|a_p|}{u_i}\right)+O\left(\sum_{u_i\leq p\leq u_{i+1}}\left|\frac{r}{u_i}-\frac{r}{p}\right|\right),\tag{3.3}$$
since $h$ is Lipschitz. By the prime number Theorem the error term above is
$$\ll \frac{r}{\log_2^2 r} \sum_{u_i\leq p\leq u_{i+1}}\frac{1}{p}\ll \frac{r}{u_i\log_2^2 r}\left(\pi(u_{i+1})-\pi(u_i)\right)\ll \frac{r}{\log r\log_2^{4} r}.$$
We assumed hypothesis D, hence we get

$$ \sum_{u_i\leq p\leq u_{i+1}}h\left(\frac{r|a_p|}{u_i}\right)= \left(\pi(u_{i+1})-\pi(u_i)\right)\int_0^Uh\left(\frac{rt}{u_i}\right)\psi(t)dt+o\left(\frac{u_i}{\log^2 u_i}\right).\tag{3.4}$$
Using that $h$ is Lipschitz again gives us
$$ \left(\pi(u_{i+1})-\pi(u_i)\right)\int_0^Uh\left(\frac{rt}{u_i}\right)\psi(t)dt= \sum_{u_i\leq p\leq u_{i+1}}\int_0^Uh\left(\frac{rt}{p}\right)\psi(t)dt+O\left(\frac{r}{\log r\log_2^{4} r}\right).$$
Therefore using the equation above with (3.3) and (3.4) give
$$ \sum_{u_i\leq p\leq u_{i+1}}h\left(\frac{r|a_p|}{p}\right)=\sum_{u_i\leq p\leq u_{i+1}}\int_0^Uh\left(\frac{rt}{p}\right)\psi(t)dt+O\left(\frac{r}{\log r\log_2^{4} r}\right).$$
Now adding all these together for $1\leq i\leq I$, we deduce
$$ \sum_{r_1\leq p\leq r_2}h\left(\frac{r|a_p|}{p}\right)= \sum_{r_1\leq p\leq r_2}\int_0^Uh\left(\frac{rt}{p}\right)\psi(t)dt+O\left(\frac{r}{\log r\log_2 r}\right).\tag{3.5}$$
For a positive number $x$ define
$$ H(x):=\int_{0}^U h\left(xt\right)\psi(t)dt, \text{ and } F(x):=\int_{0}^U f\left(xt\right)\psi(t)dt.$$
Then $$ H(x)= F(x)+ x\int_{1/x}^U t\psi(t)dt.$$
By a simple change of variables and Lemma 1.1, we have
$$ F(x)\ll\int_{0}^U f\left(xt\right)dt = \frac{1}{x}\int_{0}^{Ux}f(t)dt=\left\{\aligned &O\left(x^2\right)\ \ \ \ \ \text{ if } 0\leq x<1/U\\  &O(\log x +1) \ \text{ if } 1/U\leq x. \endaligned \right.$$
Therefore by the same analysis as in equations (1.6)-(1.9), we deduce that $$ \sum_{r_1\leq p\leq r_2}F\left(\frac{r}{p}\right)=\frac{r}{\log r}\int_0^{\infty}\frac{F(t)}{t^2}dt+o\left(\frac{r}{\log r}\right).\tag{3.6}$$
Now $$\int_0^{\infty}\frac{F(t)}{t^2}dt=\int_0^{\infty}\frac{1}{t^2}\int_{0}^U f\left(tx\right)\psi(x)dxdt=\int_0^{U}\psi(x)\int_{0}^{\infty}
\frac{f\left(tx\right)}{t^2} dtdx.$$
By a simple change of variable $T=tx$ ($t$ is the variable and $x$ is constant), we get
$$ \int_0^{\infty}\frac{F(t)}{t^2}dt=\int_0^{U}x\psi(x)\int_{0}^{\infty}
\frac{f\left(T\right)}{T^2} dTdx= N(C_2-1),\tag{3.7}$$
Furthermore
$$ \sum_{r_1\leq p\leq rU}\frac{r}{p}\int_{p/r}^Ut\psi(t)dt=N\sum_{r_1\leq p\leq rU}\frac{r}{p}-\sum_{r_1\leq p\leq rU}\frac{r}{p}\int_0^{p/r}t\psi(t)dt.\tag{3.8}$$
Now to calculate the last sum on the RHS of (3.8), we change the order of summation and integration to get
$$\sum_{r_1\leq p\leq rU}\frac{r}{p}\int_0^{p/r}t\psi(t)dt=r\int_0^{r_1/r}\sum_{r_1\leq p\leq rU}\frac{1}{p}t\psi(t)dt+r\int_{r_1/r}^U\sum_{rt\leq p\leq rU}\frac{1}{p}t\psi(t)dt.\tag{3.9}$$
And since
$$ r\int_0^{r_1/r}\sum_{r_1\leq p\leq rU}\frac{1}{p}t\psi(t)dt\ll r\log_2 r\int_0^{1/\log^2r}tdt=O\left(\frac{r}{\log^3 r}\right),\tag{3.10}$$
the LHS of (3.9) equals
$$ r\int_{r_1/r}^U\log\left(\frac{\log rU}{\log rt}\right)t\psi(t)dt +O\left(\frac{r}{\log^2 r}\right).\tag{3.11}$$
In the range $r_1/r\leq t\leq U$, we have
$$ \log\left(\frac{\log rU}{\log rt}\right)=-\log\left(1+\frac{\log t/U}{\log rU}\right)=-\frac{\log t/U}{\log rU}+O\left(\frac{\log^2 t}{\log^2 r}\right).$$
Therefore by (3.10) and (3.11) we deduce that the LHS of (3.9) equals
$$ -\frac{r}{\log rU}\int_0^Ut\psi(t)\log t/Udt+O\left(\frac{r}{\log^2 r}\right).$$
 Hence by (3.8) we have
$$
\align\sum_{r_1\leq p\leq rU}\frac{r}{p}\int_{p/r}^Ut\psi(t)dt
&=Nr\sum_{r_1\leq p\leq rU}\frac{1}{p}+ (M-N\log U)\frac{r}{\log r}+O\left(\frac{r}{\log^2 r}\right)\\
&=Nr\sum_{r_1\leq p\leq r}\frac{1}{p}+ M\frac{r}{\log r}+O\left(\frac{r}{\log^2 r}\right).\tag{3.12}\endalign$$
Thus by (3.5), (3.6), (3.7) and (3.12) we deduce that
$$ \sum_{r_1\leq p\leq r_2}h\left(\frac{r|a_p|}{p}\right)=Nr\sum_{r_1\leq p\leq r}\frac{1}{p}+(M+N(C_2-1))\frac{r}{\log r}+o\left(\frac{r}{\log r}\right).\tag{3.13}$$
Furthermore using hypothesis D and partial summation, we have that
$$\sum_{p\leq y}\frac{|a_p|-N}{p}=o\left( \frac{1}{\log y}\right),$$
and therefore we get that
$$ r\sum_{p\leq  r_1}\log |L_p(\phi_p)|=-\sum_{p\leq  r_1}rN\log \left(1-\frac{1}{p}\right)+r\log b_u+o\left(\frac{r}{\log r}\right).\tag{3.14}$$
Finally by $$\sum_{r_1<p<r}\frac{1}{p}+\log \left(1-\frac{1}{p}\right)\ll \frac{1}{r_1\log r_1},$$
and (1.10), we get the result.

\enddemo

Now define $$\Phi_u(\tau):=\pr\left(|L(1,X_u)|>b_u(e^{\gamma}\tau)^N\right).$$
Following the same lines as the proof of Theorem 1 (taking $\tau=\log( rN)+ (C_2+M/N-\log N)$) we obtain the following result

\proclaim {Theorem 3.3} With the same notations and assumptions as in Proposition 3.2, we have  for $\tau\gg 1$
$$\Phi_u(\tau)=\exp\left(-\frac{e^{\tau-C_2-M/N+\log N}}{\tau}\left(1+
o(1)\right)\right).$$

\endproclaim

Let $\{X_c(p)\}_{p\in \Cal{P}}$ be independent random variables taking the values $-1$ and $1$ with equal probabilities, and define the random Euler product
$$ L(1,X_c):=\prod_{p\in \Cal{P}}\prod_{j=1}^{d}\left(1-\frac{X_c(p)a_j(p)}{p}\right)^{-1},$$
where the $|a_j(p)|\leq 1$ are real and satisfy Hypothesis D with distribution function $\psi$. As above define $$\Phi_c(\tau):=\pr\left(L(1,X_c)>b_c(e^{\gamma}\tau)^N\right),$$ where
$$ b_{c}:= \prod_{p\in \Cal{P}}L(p)\left(1-\frac{1}{p}\right)^N  \hbox{ and } L(p)=\max_{\delta\in\{-1,1\}}\prod_{i=1}^d\left(1-\delta\frac{a_j(p)}
{p}\right)^{-1}.$$ Then using exactly the same analysis as above we prove

 \proclaim{ Theorem  3.4} Assume that the $\{a_p\}_{p\in \Cal P}$ are real and satisfy Hypothesis D. Let $M$ and $N$ be as in Theorem 4. For $r>0$ large enough we have
$$ \log\ex\left(L(1,X_c)^r\right)=rN\log_2(r)+r(\gamma N+\log b_c)+ \frac{r}{\log r}(M+N(C_1-1))+o\left(\frac{r}{\log r}\right),$$
and $$ \Phi_c(\tau)=\exp\left(-\frac{e^{\tau-C_1-M/N+\log N}}{\tau}\left(1+
o(1)\right)\right).$$

\endproclaim

\head 4. Functions of the Selberg class in the  height aspect \endhead

We begin this section by proving that the class $S^a$ contains all elements of the Selberg class which have a polynomial Euler product (all know $L$-functions do). This is a classical result which is proved in Lemmas 2.2 and 2.3 of [29].
\proclaim{Lemma 4.1}
Let $F$ be a function which satisfies axioms 1, 3 and 4, and which have a polynomial Euler product of degree $m$. Then the Ramanujan hypothesis (axiom 2) is equivalent to the fact that the local roots of $F$ satisfy $|\alpha_{i,F}(p)|\leq 1$, for all prime $p$ and $1\leq i\leq m$. Moreover we have that $S^p\subset S^a$.
\endproclaim
\demo{Proof} By the identity between the Dirichlet series and the Euler product of $F$ (which holds in Re$(s)>1$), we have the power series expansion
$$ \sum_{k=1}^{\infty}a_F(p^k)X^k
=\prod_{i=1}^m\left(1-\alpha_{i,F}(p)X\right)^{-1}.$$
Therefore if the Ramanujan hypothesis (axiom 2) is true, then the series is analytic for $|X|<1$, hence so is the right-hand side which implies that $|\alpha_{i,F}(p)|\leq 1$.

Now if $|\alpha_{i,F}(p)|\leq 1$ for all prime $p$ then one can easily show that $|a_F(n)|\leq d_m(n)\ll n^{\epsilon}$ for all $\epsilon$, where $d_m(n)$ is the $m$-th divisor function. Furthermore in this case, we have
$$\sum_{n\leq x}|a_F(n)|\leq \sum_{n\leq x}d_m(n)\ll x(\log x)^{m-1},$$ which implies that $S^p\subset S^a.$

\enddemo
Let $F$ be a function which satisfies axioms 1, 3 and 4. In order to study values of $F$ on the line Re$(s)=1$, we will need good bounds for $\sum_{p\leq x}|a_F(p)|/p$. If $F\in S^p$ and has degree $m$, then it is easy to see that the last sum is bounded by $(m+o(1))\log\log x$, which follows from the fact that $|a_F(p)|\leq m$. While if $F\in S^a$, then by partial summation we can show that
 $$ \sum_{p\leq x}\frac{|a_F(p)|}{p}\ll (\log x)^{\beta_F+1}.\tag{4.1}$$

 To study the distribution of large values of $|F(1+it)|$, the first step is to compute its moments. However it is difficult in general to find asymptotic formulas for large moments due to the difficulty of estimating the off-diagonal terms. To solve this problem we approximate $F(1+it)$ by very short Euler products and compute large moments of these. We generalize the approach of Granville and Soundararajan ([7] and [8]), and use assumption (WZD) to show that this approximation holds for almost all $t$. First we prove

\proclaim { Lemma 4.2} Let $F\in S^a$, and $y\ge 2$ and $|t| \ge y+3$ be real numbers. Let $1/2\leq \sigma_0<\sigma$, and suppose that there are no
zeros of $F(s)$ inside the rectangle $\{z: \sigma_0\leq
\text{Re}(z)\leq 1, |\text{Im}(z)-t|\leq y+1\}$, then for
$s=\sigma+it$ we have
$$|\log F(s)|\ll \frac{\log|t|}{\sigma-\sigma_0}$$
\endproclaim

\demo{Proof} First for $\sigma> 1$ we know that $|\log F(s)|\ll 1$
so there is nothing to prove in this case. Now assume that
$\sigma\leq 1$, and consider the circles with center $2+it$ and
radii $r:=2-\sigma<R:=2-\sigma_0$, so that the smaller circle passes
through $s$. By our hypothesis $\log F(s)$ is analytic inside the
larger circle. For a point $z$ on the larger circle we use the
bound $|F(z)|\ll |z|^{\alpha_F} \ll |t|^{\alpha_F}$, which holds for some constant $\alpha_F$ by the functional equation, so that $\text{Re}(\log F(z))\ll \log|t|$. Thus upon using the
Borel-Caratheodory Theorem we deduce that
$$ |\log F(s)|\leq \frac{2r}{R-r}\max_{|z-2-it|=R} \text{Re} (\log
F(s)) +\frac{R+r}{R-r}|\log
F(2+it)|\ll\frac{\log|t|}{\sigma-\sigma_0}.$$
\enddemo
Using this bound we prove
\proclaim {Lemma 4.3}  Let $F\in S^a$, and $y\ge 2$ and $|t| \ge y+3$ be real numbers.
Let $\frac{1}{2} \leq \sigma_0 <\sigma\leq 1$ and suppose that the
rectangle $\{ s: \ \ \sigma_0 <\text{Re}(s) \leq 1, \ \
|\text{Im}(s) -t| \leq y+2\}$ does not contains zeros of $F(s)$.
Then
$$
\log F(\sigma+it)= \sum_{n=2}^{y} \frac{b_n}{n^{\sigma+it}} +
O\left( \frac{\log
|t|}{(\sigma_1-\sigma_0)^2}y^{\sigma_1-\sigma}\right),
$$
where  $\sigma_1 = \min(\sigma_0+\frac{1}{\log y},
\frac{\sigma+\sigma_0}{2})$.
\endproclaim
\demo{Proof} Without loss of generality we may assume that $y\in
{\Bbb Z}+\frac{1}{2}$. Let $c= 1-\sigma+\frac{1}{\log y}$, by
Perron's formula (see Davenport [2]), we have
$$
\align \frac{1}{2\pi i} \int_{c-iy}^{c+iy} \log
F(\sigma+it+w)\frac{y^w}{w} dw &= \sum_{n=2}^{y}
\frac{b_n}{n^{\sigma+it}} + O\Big(\frac{1}{y} \sum_{n=1}^{\infty}
\frac{y^c|b_n|}{n^{\sigma+c}} \frac{1}{|\log (y/n)|}\Big)\tag{4.2}
\endalign
$$
To deal with the error term in the RHS of (4.2), we take first all the terms
 for which $n\leq \frac{3}{4}y$ or $n\geq \frac54 y$. For these, $|\log (y/n)|$
  has a positive lower bound, so their contribution is
  $$\ll \frac{1}{y^{\sigma}} \sum_{n=1}^{\infty}\frac{|b_n|}{n^{1+1/\log y}}=
   \frac{1}{y^{\sigma}}\left(\sum_{p \in {\Cal P}}\frac{|a_F(p)|}{p^{1+1/\log y}}+O(1)\right)\ll \frac{\log y}{y^{\sigma}},$$
by (4.1).

Now consider the terms for which $\frac34 y<n<y$. Let $y_0$
be the largest prime power less than $y$; we can suppose that
$\frac34 y<y_0<y$. For the term $n=y_0$, we have that $\log(y/n)\geq (y-n)/y$, so the contribution of this term is $\ll \frac{|b_n|}{y^{\sigma}}\ll y^{\lambda-\sigma}$ by axiom 4. For the other terms we can put $n=y_0-m$, where $0<m<y/4$, and then $\log(y/n)\geq m/y$. Hence the contribution of these terms is
$$\ll \frac{1}{y^{\sigma}}\sum_{0<m<y/4}\frac{|b_{y_0-m}|y_0}{(y_0-m)^{1+1/\log y}m}\ll y^{\lambda-\sigma}\log y.$$
The terms with $ y<n<\frac54y$ are dealt with similarly, except that $y_0$ is replaced by $y_1$ which is the least prime power greater than $y$. Thus the error term in the RHS of (4.2) is
$$\ll y^{\lambda-\sigma}\log y.$$

 Now we move the line of integration to the left at Re$(w)=\sigma_1
-\sigma <0$.  From our hypothesis we encounter only a simple pole at
$w=0$ which has residue $\log F(\sigma+it)$. Thus the left side of
(4.2) equals $\log F(\sigma+it)$ plus
$$
\frac{1}{2\pi i} \Big( \int_{c-iy}^{\sigma_1-\sigma-iy}
+\int_{\sigma_1-\sigma -iy}^{\sigma_1-\sigma+iy}
+\int_{\sigma_1-\sigma +iy}^{c+iy}\Big) \log
F(\sigma+it+w)\frac{y^w}{w} dw \ll \frac{\log
|t|}{(\sigma_1-\sigma_0)^2}y^{\sigma_1-\sigma},
$$
by Lemma 4.2.

\enddemo
\proclaim{Lemma 4.4} Let $F\in S^a$ and satisfies assumption (WZD). Let $T>0$ be large. Then there exists some $\nu_F>0$, such that for any
real number $3\leq y\leq T/2$ the asymptotic
$$\log F(1+it)=\sum_{n=2}^y\frac{b_n}{n^{1+it}}
+O\left(y^{-\nu_F}\log^3 T\right)$$ holds for all $t\in (T,2T)$ except
a set of measure $T^{1-c_F}y$, for some fixed constant $c_F>0$.
\endproclaim

\demo{Proof} By (WZD), there exist two positive constants $0<\nu_F,c_F<1/2$ such that

\noindent $N_F(1-\nu_F,T)\ll T^{1-c_F}$.  Thus appealing  to Lemma 4.3 with $\sigma_0= 1-\nu_F$ gives the result.
\enddemo

Let $F\in S^a$ and satisfies (WZD). For $y>2$ define the short Euler product of length $y$ of $F(1+it)$ by
$$F(1+it,y):=\sum_{n\in S(y)}\frac{a_F(n) }{n^{1+it}} , \text{ where } S(y)=\{n\in {\Bbb N}: p|n\implies p\leq y\}.\tag{4.3}$$
And similarly define the short Euler random product of $F(1,X)$ by
$$F(1,X,y):=\sum_{n\in S(y)}\frac{a_F(n) X(n)}{n}.$$
The following proposition establish that large integral moments of $F(1+it,y)$ and $F(1,X,y)$ are roughly equal. For simplicity we let $a_F(n)=a_n$.

\proclaim{ Proposition 4.5} Let $F\in S^a$.  Let $T>0$ be large and $2<y<(\log
T)^A$ a real number where $A>0$ is a fixed constant. Then for
all positive integers $k$ in the range $1\leq k\leq \log
T/(B(\log_2T)^{\beta_F+2})$ (for some suitably large constant $B=B(F,A)>0$), we have
$$\frac{1}{T}\int_T^{2T}|F(1+it,y)|^{2k}dt = {\Bbb E}\left(|F(1,X,y)|^{2k}\right)+o(1).$$
Moreover if $F\in S^p$, then this last estimate holds in the wider range

\noindent $1\leq k\leq \log
T/(C\log_2T\log_3 T)$ for some suitably large constant $C=C(F,A)>0$.
\endproclaim
\demo {Proof} First we have that

$$|F(1+it,y)|^{2k}=\sum_{n,m\in S(y)} \frac{c_k(n)\overline{c_k(m)}}{m^{1-it}n^{1+it}},\ \text{ where } \ c_k(m)=\sum_{m_1m_2...m_k=m}a_{m_1}...a_{m_k}.$$
Therefore
$$\frac{1}{T}\int_T^{2T}|F(1+it,y)|^{2k}dt=\sum_{n,m\in S(y)}
\frac{c_k(n)\overline{c_k(m)}}{mn}\frac{1}{T}\int_T^{2T}\left(\frac{m}{n}\right)^{it}dt.$$
The contribution of the diagonal terms $m=n$ is
$$ \sum_{n\in S(y)} \frac{|c_k(n)|^{2k}}{n^2}={\Bbb E}\left(|F(1,X,y)|^{2k}\right),$$
so it remains only to bound the contribution of the non-diagonal
terms. To do so we divide the above sum into two parts: $E_1$ over
the integers $1\leq m\neq n\leq \sqrt{T}$, and $E_2$ over the
remaining terms. For the first terms we use the fact that
$$\int_T^{2T}\left(\frac{m}{n}\right)^{it}dt\ll
\frac{1}{T|\log(m/n)|}\leq \frac{1}{\sqrt{T}},$$ which implies that
$$ E_1\ll \frac{1}{\sqrt{T}}\left(\sum_{n\in
S(y)}\frac{|c_k(n)|}{n}\right)^2=\frac{1}{\sqrt{T}}\left(\sum_{n\in
S(y)}\frac{|a_n|}{n}\right)^{2k}= \frac{1}{\sqrt{T}}\exp\left(2k\sum_{p\leq y}\frac{|a_p|}{p}+ O(k)\right).$$  Let $\beta=1/\log_2 T$. The
contribution of the remaining non-diagonal terms is
$$ E_2\ll \left(\sum\Sb n>\sqrt{T}\\ n\in S(y)\endSb\frac{|c_k(n)|}{n}\right)\left(\sum_{n\in
S(y)}\frac{|c_k(n)|}{n}\right)\ll \frac{1}{T^{\beta/2}}\left(\sum_{n\in
S(y)}\frac{|c_k(n)|}{n^{1-\beta}}\right)\left(\sum_{n\in
S(y)}\frac{|c_k(n)|}{n}\right).$$ Moreover
$$ \left(\sum_{n\in
S(y)}\frac{|c_k(n)|}{n^{1-\beta}}\right)=\exp\left(k\sum_{p\leq
y}\frac{|a_p|}{p^{1-\beta}}+O(k)\right)=\exp\left(O\left(k\sum_{p\leq
y}\frac{|a_p|}{p}\right)\right).$$ Therefore the contribution of the non-diagonal terms
is bounded by $$\exp
\left(-\beta\log T+ O\left(k\sum_{p\leq
y}\frac{|a_p|}{p}\right)\right).$$
Finally the result follows from the fact that
$$ \sum_{p\leq
y}\frac{|a_p|}{p}\ll \left\{\aligned & \log\log y \ \ \ \text{ if } F\in S^p,\\
& (\log y)^{\beta_F+1} \ \ \ \text{ if } F\in S^a.\\
\endaligned \right.$$
\enddemo

By Lemma 4.4 we can approximate $F(1+it)$ by very short Euler products for almost all $t\in[T,2T]$, where the exceptional set have measure $\leq T^{1-c_F}$ for some positive constant $c_F$. Therefore in order to compute high moments of $F(1+it)$, we need to control the contribution of these exceptional values. To do so we prove the classical bound for $F(1+it)$. We follow closely the method of  Soundararajan [28].
\proclaim {Lemma 4.6} Let $F\in S^a$. Then there exists a constant $m_F$ such that for any real number $t$ we have
$$ F(1+it)\ll \log (|t|+2)^{m_F}.$$
\endproclaim

\demo{Proof} Let $c>0$, and take $Z\geq 1$. Then expand $F(1+it+\omega)$ into its Dirichlet series to see that
$$ \frac{1}{2\pi i}\int_{c-i\infty}^{c+i\infty}F(1+it+\omega)Z^{\omega}\Gamma(\omega)d\omega= \sum_{n=1}^{\infty}\frac{a_n}{n^{1+it}}e^{-n/Z}.$$
Now shift the contour to the line of integration $\omega=-1+\epsilon$. the pole at $\omega =0$ leaves the residue $F(1+it)$. The possible pole at $\omega =-it$ leaves a residue $\ll Z^{\epsilon}(1+|t|)^{\epsilon}|\Gamma(it)|\ll Z^{\epsilon}e^{-|t|}$. Finally using the functional equation (axiom 3) and Stirling's formula we have
$$ \left|F(1+it+\omega)\right|= \frac{|\overline{\gamma_F}(-it-\omega)|}{|\gamma_F(1+it+\omega)|}\left|F(-it-\omega)\right|\ll (1+|t|+|\omega|)^{1+\epsilon},$$
for $\omega$ on the line Re$(\omega)=-1+\epsilon.$ Hence the integral on the line Re$(\omega)=-1+\epsilon$ is $\ll Z^{-1+\epsilon}(1+|t|)^{1+\epsilon}$.
Therefore we deduced that for any $\epsilon>0$ and $Z\geq 1$
$$F(1+it)=\sum_{n=1}^{\infty}\frac{a_n}{n^{1+it}}e^{-n/Z}+O\left(Z^{-1+\epsilon}
(1+|t|)^{1+\epsilon}+Z^{\epsilon}e^{-|t|}\right).$$
Choose $Z=1+|t|^2$. Then the error term is $\ll |t|^{-1+2\epsilon}$. Now the main term is bounded by
$$ \sum_{n\leq Z^2}\frac{|a_n|}{n}+\sum_{k=1}^{\infty}e^{-Z^{k}}\sum_{Z^{k+1}\leq n\leq Z^{k+2}} \frac{|a_n|}{n}.$$
Finally since $F\in S^a$, we know by (3) that $\sum_{n\leq y}|a_n|/n\leq (\log y)^{\beta_F+1}$. Therefore we get
$$ F(1+it)\ll \log (|t|+2)^{m_F},$$ for some constant $m_F>0$.

\enddemo

\demo{Proof of Theorem 3}
 Let $\nu_F$ be as in Lemma 4.4 and take $y=\log^{(A+3)/\nu_F} T.$ Then by this Lemma we know that
$$ F(1+it)= F(1+it,y)\left(1
+O\left(\frac{1}{\log^A T}\right)\right),$$
for all $t\in [T,2T]$ except a set of measure $V=T^{1-\delta}$ for some small $\delta=\delta(F,A).$
Now by Lemma 4.6 we deduce that
$$\frac{1}{T}\int_T^{2T}|F(1+it)|^{2k}dt= \frac{1}{T}\int_T^{2T}|F(1+it,y)|^{2k}dt\left(1
+O\left(\frac{1}{\log^A T}\right)\right)+ Er,$$
where $$ Er\ll V\max_{t\in [T,2T]}|F(1+it)|^{2k}+V \max_{t\in [T,2T]}|F(1+it,y)|^{2k}\ll T^{1-\delta}\exp\left(O(k \log_2 T)\right)=o(1).$$
Furthermore by (3.2) and the fact that $h(t)\ll t^2$ for $0\leq t\leq 1$, we have
$$
\align
{\Bbb E}\left(|F(1,X)|^{2k}\right)&= {\Bbb E}\left(|F(1,X,y)|^{2k}\right)\exp\left(O\left(k\sum_{p>y}\frac{|a_p|^2}{p^2}\right)\right)\\
&= {\Bbb E}\left(|F(1,X,y)|^{2k}\right)\left(1
+O\left(\frac{1}{\log^A T}\right)\right).\\
\endalign
$$
Therefore by Proposition 4.5 we get the result.
\enddemo
We now deduce information on large values of $|F(1+it)|$. Indeed we prove Corollary 1
\demo{Proof of Corollary 1}
By Theorem 3 and Lemma 3.1 we have that
$$  \frac{1}{T}\int_T^{2T}|F(1+it)|^{2k}dt=\exp\left(2k\sum_{p\leq 2k}\frac{|a_p|}{p}+O(k)\right),$$
uniformly for $k\leq (\log T)/B(\log_2 T)^{\beta_F+2}$. Therefore taking $k=[(\log T)/B(\log_2 T)^{\beta_F+2}]$, we have
$$\max_{t\in [T,2T]}|F(1+it)|\geq \exp\left(\sum_{p\leq 2k}\frac{|a_p|}{p}+O(1)\right),$$
from which we deduce the result.
\enddemo
In order to prove Theorem 4 we consider only those $F\in S^p$ which satisfy Hypothesis D (with distribution function $\psi$). Hence our random model $F(1,X)$ satisfies the assumptions of Proposition 3.2, and thus we can compute its real moments. Using this fact along with the same saddle point method used to prove Theorem 1 we prove Theorem 4.

\demo{Proof of Theorem 4} First for any positive real number $r$ we have that
$$ rN \int_0^{\infty}\Phi_F(t) t^{rN-1}dt =\left(b_Fe^{\gamma N}\right)^{-r}\frac{1}{T}\int_T^{2T}|F(1+it)|^{r}dt.$$
Take $r=2k$ where $1\leq k\leq \log T/(B\log_2 T\log_3 T)$ is a positive integer, and $B$ is a suitably large constant.  Therefore by Theorem 3 and Proposition 3.2 we have that
$$
\align
2kN \int_0^{\infty}\Phi_F(t) t^{2kN-1}dt &=\left(b_Fe^{\gamma N}\right)^{-2k}{\Bbb E}\left(|F(1,X)|^{2k}\right)\left(1+O\left(\frac{1}{\log T}\right)\right)\\
&= (\log 2k)^{2kN}\exp\left(\frac{2k}{\log 2k}(M+N(C_2-1))+o\left(\frac{2k}{\log 2k}\right)\right).\tag{4.4}\\
\endalign$$
In order to use the saddle point method as in the proof of Theorem 1, we require (4.4) to be valid for all real numbers $r$ in the range  $2\leq r\leq 2\log T/(B\log_2 T\log_3 T)$. To do so we know by (4.4) (for $k\geq 1/N$ and fixed) that
$$ \int_0^{\infty} \Phi_F(t)dt\ll 1+\int_1^{\infty} \Phi_F(t)t^{2kN-1}dt\ll 1.$$
Therefore by H\"older's inequality we have for $a<b$
$$ \int_0^{\infty} \Phi_F(t)t^adt\leq \left(\int_0^{\infty} \Phi_F(t)dt)^{1-a/b}\right)\left(\int_0^{\infty} \Phi_F(t)t^bdt\right)^{a/b}\ll\left(\int_0^{\infty} \Phi_F(t)t^bdt\right)^{a/b} .$$
Using this we may interpolate (4.4) to non-integer value $\kappa\in (k-1,k)$ by taking $a=2(k-1)N-1, b=2\kappa N-1$ and then $a=2\kappa N-1, b=2kN-1$ in the last inequality to obtain
$$\left(\int_0^{\infty} \Phi_F(t)t^{2(k-1) N-1}dt\right)^{\frac{2\kappa N-1}{2(k -1)N-1}}\ll \int_0^{\infty} \Phi_F(t)t^{2\kappa N-1}dt\ll \left(\int_0^{\infty} \Phi_F(t)t^{2kN-1}dt\right)^{\frac{2\kappa N-1}{2k N-1}}.$$
Hence we get (4.4) for $r=2\kappa$ by substituting (4.4) for $2(k-1)$ and $2k$ into this equation. Finally applying the same saddle point method as in the proof of Theorem 1 (now taking $\tau=\log r+C_2+M/N$) we get the result.

\enddemo

\head 5. Distribution of Quadratic twists of automorphic $L$-functions \endhead

Let $\pi$ be an auto-dual cuspidal automorphic representation of $GL(m)/{\Bbb Q}$, and let $L(\pi, \cdot)$ be the $L$-function attached to $\pi$. Consider the values of  $L(\pi\otimes \chi_d,1)$ as $d$ varies over all fundamental discriminants with $|d|\leq x$, where $\chi_d$ denotes the unique real character mod $d$. The probabilistic model constructed for this family is based on the totally multiplicative function $X(.)$, where the $\{X(p)\}_{p \in \Cal{P}}$ are independent random variables taking the values $1$ with probability $p/(2(p+1))$, $0$ with probability $1/(p+1)$, and $-1$ with probability $p/(2(p+1))$; and for a positive integer $n=p_1^{a_1}...p_r^{a_r}$ we have $X(n):=X(p_1)^{a_1}...X(p_r)^{a_r}.$ Then define the following random series
$$ L_{\pi}(1,X):=\sum_{n=1}^{\infty}\frac{a_{\pi}(n)X(n)}{n},  \hbox{ (these series converge with probability }1).$$
As in equation (4.3) let $L(\pi\otimes \chi_d,1,y)$ and $L_{\pi}(1,X,y)$ denote the short Euler products of length $y$ of $L(\pi\otimes \chi_d,1)$ and $L_{\pi}(1,X)$ respectively.

The analogue of (WZD) (zeros density estimates near Re$(s)=1$) has been proved by Kowalski and Michel (see Theorem 5 of [14]) for general families $\Cal{F}$ of cuspidal automorphic representations of $GL(m)/{\Bbb Q}$, under the condition that $|a_{f}(p)|\ll p^{\theta}$ for some $0<\theta<1/4$ and all $f \in \Cal{F}.$  For our family this holds if $$|a_{\pi}(p)|\ll p^{\theta} \text{ for some } 0<\theta<1/4. \tag{5.1}$$ In particular this is known to be true for $GL(2)$ Maass cusp forms, where one can take $\theta=1/9$ by the work of Kim and Shahidi [13].

The analysis for this section is similar to section 4, indeed the methods for proving the analogues of Lemmas 4.2, 4.3, 4.4 and 4.6 are exactly the same. However the method of proving the estimates for the moments of short Euler products is a little bit different so we include it.

\proclaim{ Proposition 5.1} Let $\pi$ be an auto-dual cuspidal automorphic representation of $GL(m)/{\Bbb Q}$. Let $x>0$ be large and $2<y<(\log
x)^A$ a real number where $A>0$ is a fixed constant. If $\pi$ satisfies assumption (3), then for
all positive integers $k$ in the range $1\leq k\leq \log
x/(B_1(\log_2x)^{\beta_{\pi}+2})$ (for a suitably large constant $B_1=B_1(\pi,A)$), we have

$$\frac{\pi^2}{6x}\sums_{|d|\leq x}L(\pi\otimes \chi_d,1,y)^{k} = {\Bbb E}\left(L_{\pi}(1,X,y)^{k}\right)+o(1),$$
where $\sums$ denotes the sum over fundamental discriminants. Moreover if $\pi$ satisfies the Ramanujan hypothesis then this last estimate holds in the wider range

\noindent $1\leq k\leq \log
x/(B_2\log_2x\log_3 x)$ (for a suitably large constant $B_2=B_2(\pi,A)$).
\endproclaim

\demo{Proof}
For simplicity we assume that $\pi$ satisfies the Ramanujan hypothesis, since the proof is similar in the other case (where $\pi$ satisfies (3)), exactly as in the proof of Proposition 4.5. In this case we have $|a_{\pi}(p)|\leq m$ by Lemma 4.1. First we have
$$ L(\pi\otimes \chi_d,1,y)^{k}=\left(\sum_{n\in S(y)}\frac{a_{\pi}(n)\chi_d(n)}{n}\right)^k=\sum_{n\in S(y)}\frac{c(k,n)\chi_d(n)}{n},$$
where $c(k,n)=\sum_{n_1n_2...n_k=n}a_{\pi}(n_1)...a_{\pi}(n_k)$. Further one can see that $\ex(X(n))=\prod_{p|n}\left(\frac{p}{p+1}\right)$ if $n$ is a square and equals $0$ otherwise. Hence
$$ \ex\left(L_{\pi}(1,X,y)^k\right)=\ex\left(\sum_{n\in S(y)}\frac{c(k,n)X(n)}{n}\right) = \sum_{l\in S(y)}\frac{c(k,l^2)}{l^2}\prod_{p|l}\left(\frac{p}{p+1}\right)
.$$
Now we estimate the contribution of the diagonal terms $n=\square$ (which constitute the main term) to the moments of $L(\pi\otimes \chi_d,1,y)$. When $n=l^2$ we have
$$ \sums_{|d|\leq x}\chi_d(l^2)=\sums_{\Sb |d|\leq x\\ (l,d)=1\endSb}1=\frac{6}{\pi^2}x\prod_{p|l}\left(\frac{p}{p+1}\right)+O\left(
x^{1/2+\epsilon}d(l)\right).$$ Thus the contribution of the diagonal terms is
$$\frac{6}{\pi^2}x\sum_{l\in S(y)}\frac{c(k,l^2)}{l^2}\prod_{p|l}\left(\frac{p}{p+1}\right)
+O\left(x^{1/2+\epsilon}\sum_{l\in S(y)}\frac{c(k,l^2)d(l)}{l^2}\right). \tag{5.2} $$
Before considering the non-diagonal terms, we will bound the error term in the RHS of (5.2). To do so we divide the sum into two parts $l\leq x^{1/3}$ and $l>x^{1/3}$. Using that $d(l)\leq \sqrt{l}$, the contribution of the first terms is
$$
\align &\ll x^{2/3+\epsilon}\sum_{l\in S(y)}\frac{c(k,l^2)}{l^2}\leq x^{2/3+\epsilon}\left(\sum_{n\in S(y)}\frac{|a_{\pi}(n)|}{n}\right)^k\\
&\leq x^{2/3+\epsilon}\exp\left(k\sum_{p\leq y}\frac{|a_{\pi}(p)|}{p}+O(k)\right) \ll x^{3/4},\tag{5.3}
\endalign
$$ by our hypothesis on $k$. Now for the other terms we use the bound of Ramanujan
$$d(l)\leq \exp\left(\frac{\log 2\log(l)}{\log_2(l)}(1+o(1))\right)\leq l^{2\beta}, $$
with $\beta=1/(2\log_2 x)$. Thus similarly to (5.3), we can show that the contribution of these terms is
$$
\ll x^{1/2+\epsilon}\sum_{l\in S(y)}\frac{c(k,l^2)}{l^{2-2\beta}}
\leq x^{1/2+\epsilon}\exp\left(k\sum_{p\leq y}\frac{|a_{\pi}(p)|}{p^{1-\beta}}+O(k)\right)\ll x^{2/3} ,\tag{5.4}
$$
 using that $p^{\beta}=O(1)$ for all primes $p\leq y$. Moreover the contribution of the non-diagonal terms $n\neq \square$ is

$$\ll \sum_{\Sb n\in S(y)\\n\neq \square\endSb}\frac{|c(k,n)|}{n}\bigg|\sums_{|d|\leq x}\chi_d(n)\bigg|.$$
Consider first the terms $n\leq x$. For these we use Lemma 4.1 of Granville-Soundararajan [8] which asserts that
$$\bigg|\sums_{|d|\leq x}\chi_d(n)\bigg|\ll x^{1/2}n^{1/4}(\log n)^{1/2},$$
if $n$ is not a perfect square.
Hence the contribution of these terms is
$$ \ll x^{5/6}\sum_{n\in S(y)}\frac{|c(k,n)|}{n}\ll x^{5/6}\exp\left(k\sum_{p\leq y}\frac{|a_{\pi}(p)|}{p}+O(k)\right)\ll x^{6/7}.\tag{5.5}$$
Now the contribution of the terms $n>x$ is
$$ \ll x \sum_{\Sb n\in S(y)\\ n>x \endSb}\frac{|c(k,n)|}{n}\ll x^{1-2\beta} \sum_{n\in S(y)}\frac{|c(k,n)|}{n^{1-2\beta}}\ll x^{1-2\beta}\exp\left(O\left(k\sum_{p\leq y}\frac{|a_{\pi}(p)|}{p}\right)\right)\ll x^{1-\beta}.$$
This along with (5.2)-(5.5) complete the proof.
\enddemo
Now using exactly the same method as in the proof of Theorem 3, we can prove the following result
\proclaim {Theorem 5.2}
Let $\pi$ be an auto-dual cuspidal automorphic representation of $GL(m)/{\Bbb Q}$. Let $x>0$ be large, and take $A>0$. If $\pi$ satisfies assumptions (3) and (5.1), then for all positive integers $k$ in the range $1\leq k\leq \log
x/(B_1(\log_2x)^{\beta_{\pi}+2})$ (for a suitably large constant $B_1=B_1(\pi,A)$), we have

$$\frac{\pi^2}{6x}\sums_{|d|\leq x}L(\pi\otimes \chi_d,1)^{k} = {\Bbb E}\left(L_{\pi}(1,X)^{k}\right)\left(1+O\left(\frac{1}{\log^A x}\right)\right).$$
Moreover if $\pi$ satisfies the Ramanujan hypothesis then this last estimate holds in the wider range $1\leq k\leq \log
x/(B_2\log_2x\log_3 x)$ (for a suitably large constant $B_2=B_2(\pi,A)$).
\endproclaim
From this theorem and using the same proof as for Corollary 1, we deduce the following result on large values of $L(\pi\otimes\chi_d,1)$
\proclaim{ Corollary 5.3 }
Let $\pi$ be an auto-dual cuspidal automorphic representation of $GL(m)/{\Bbb Q}$ satisfying the Ramanujan hypothesis and assumption (8). Then for $x>0$ large, there is a fundamental discriminant $d$ with $|d|\leq x$ such that
$$L(\pi\otimes \chi_d,1)\gg\left\{\aligned &(\log\log x)^{\kappa_{\pi}+o(1)}\ \ \ \ \ \text{ if } E_{\pi}(t)=o(\log\log t),\\  &(\log\log x)^{\kappa_{\pi}} \ \text{ if } E_{\pi}(t)=O(1). \endaligned \right.$$
\endproclaim

Finally we prove Theorem 5
\demo {Proof of Theorem 5}
 We use exactly the same proof as Theorem 4. What remains only is compute the moments of $L_{\pi}(1,X)$. To this end let $\{X_1(p)\}_{p\in\Cal{P}}$ be independent random variables taking the values $-1$ and $1$ with equal probabilities. Extend $X_1(.)$ to a completely multiplicative function  over positive integers and consider the random series $$L_{\pi}(1,X_1):=\sum_{n=1}^{\infty}\frac{a_{\pi}(n)X_1(n)}{n}.$$
 Then using Theorem 3.4 we can compute the moments of $L_{\pi}(1,X_1)$. Now since ${\Bbb E}(X)={\Bbb E}(X_1)=0$ we have that
$$\align\frac{{\Bbb E}\left(L_{\pi}(1,X)^{k}\right)}{{\Bbb E}\left(L_{\pi}(1,X_1)^{k}\right)}&= \frac{{\Bbb E}\left(L_{\pi}(1,X,k^2)^{k}\right)}{{\Bbb E}\left(L_{\pi}(1,X_1,k^2)^{k}\right)}\exp\left(O\left(k\sum_{p>k^2}\frac{|a_{\pi}(p)|^2}{p^2}\right)\right)\\
&= \prod_{p\leq k^2}\left(1+O\left(\frac{1}{p}\right)\right)\left(1+O\left(\frac{1}{k}\right)\right)=\log k^{O(1)}.\endalign$$

 Thus by Theorem 3.4 the asymptotic of $\log{\Bbb E}\left(L_{\pi}(1,X_1)^{k}\right)$  holds also for  $\log{\Bbb E}\left(L_{\pi}(1,X)^{k}\right)$, completing the proof.

\enddemo

\enddocument